\documentclass[11pt]{article}
\usepackage[dvips]{graphics}
\input{amssym.def}
\input{amssym.tex}

\title{\Large\bf  
Sadovskii vortex  in a wedge and \\
the associated Riemann-Hilbert problem on a torus\footnote{The research was partly supported by EPSRC grant no EP/R014604/1 and the Simons Foundation grants no. 319217 and 713080.}}
\author{\bf {\sc Y.A. Antipov and  A.Y. Zemlyanova}\\ 
Department of Mathematics, Louisiana State University\\
Baton Rouge, LA 70803, USA\\
Department of Mathematics, Kansas State University\\
Manhattan, KS 66506,USA}

\date{}

\newcommand{\sgn}{\mathop{\rm sgn}\nolimits}

\newcommand{\I}{\mathop{\rm Im}\nolimits}
\newcommand{\R}{
\mathop{\rm Re}\nolimits}
\newcommand{\const}{\mbox{const}}

\newcommand{\ov}[1]{\overline{#1}}

\newcommand{\Ga}{\alpha}
\newcommand{\Gb}{\beta}
\newcommand{\Gd}{\delta}

\newcommand{\Gf}{\phi}

\newcommand{\Gg}{\gamma}
\newcommand{\Gc}{\chi}

\newcommand{\Gl}{\lambda}
\newcommand{\Gn}{\eta}
\newcommand{\Gm}{\mu}

\newcommand{\Gr}{\rho}

\newcommand{\Gs}{\sigma}

\newcommand{\Gj}{\tau}

\newcommand{\Go}{\omega}
\newcommand{\Gx}{\xi}

\newcommand{\Gz}{\zeta}

\newcommand{\GF}{\Phi}
\newcommand{\GG}{\Gamma}
\newcommand{\GL}{\Lambda}

\newcommand{\GO}{\Omega}

\newcommand{\CA}{{\cal A}}
\newcommand{\CB}{{\cal B}}

\newcommand{\CD}{{\cal D}}

\newcommand{\CF}{{\cal F}}

\newcommand{\CI}{{\cal I}}
\newcommand{\CJ}{{\cal J}}

\newcommand{\CL}{{\cal L}}

\newcommand{\CP}{{\cal P}}

\newcommand{\CR}{{\cal R}}

\newcommand{\CW}{{\cal W}}

\def\Ba{{\bf a}}
\def\Bb{{\bf b}}

\def\Bp{{\bf p}}
\def\Bq{{\bf q}}

\def\Bv{{\bf v}}

\def\BE{{\bf E}}

\def\BK{{\bf K}}

\newcommand{\beq}{\begin{equation}}
\newcommand{\eeq}{\end{equation}}
\newcommand{\barr}{\begin{eqnarray}}
\newcommand{\earr}{\end{eqnarray}}
\newcommand{\beqn}{\begin{equation*}}
\newcommand{\eeqn}{\end{equation*}}
\newcommand{\barrn}{\begin{eqnarray*}}
\newcommand{\earrn}{\end{eqnarray*}}

\newcommand{\fr}{\frac}

\begin{document}
\maketitle

\begin{abstract}

Reconstruction of conformal mappings from canonical slit domains onto multiply-connected physical domains with 
a free boundary is of interest in many different models arising in fluid mechanics. In the present paper, an exact formula for the conformal map from the exterior of two slits onto the doubly connected flow domain is obtained when a fluid flows in a wedge about a Sadovskii vortex. The map is employed to determine the potential
flow outside the vortex and the  vortex domain boundary
provided the circulation $\GG$ around the vortex and constant speed $U$ on the vortex boundary are prescribed, and there are no stagnation points on the walls. The map is expressed in terms of a rational function on an elliptic surface topologically equivalent to a torus and the solution to a symmetric Riemann-Hilbert problem on a finite and a semi-infinite segments on the same genus-1 Riemann surface. Owing to its special features, the Riemann-Hilbert problem requires a  novel analogue of the Cauchy kernel on an elliptic surface. Such a kernel is proposed and employed to derive a closed-form solution to the Riemann-Hilbert problem and  the associated Jacobi inversion problem. The final formula for the conformal map possesses a free geometric parameter and two model parameters, the wedge angle $\Ga$ and $\GG/U$. It is shown that when $\Ga<\pi$ the solution exists and the vortex has two cusps, while  the solution does not exist when the wedge angle exceeds $\pi$.

\end{abstract}

\setcounter{equation}{0}

\section{Introduction}

Free  boundary problems have attracted considerable attention 
because of their challenge and
numerous applications including  Hele-Show flow \cite{gus}, supercavitating flow
\cite{bir}, flow around vortex patches \cite{Saffman1992},
 inverse elastic problems  of inclusions whose stress field is prescribed \cite{esh}, hypersonic flow in air around a 
body when a realistic model requires to determine shock shapes, the boundary layer, 
the recirculation region, and the wake behind the body \cite{vag}.

In this paper we further  advance  the method of the Riemann-Hilbert problem on Riemann surfaces
for conformal maps of multiply connected free boundary  domains    \cite{AntipovSilvestrov2007, az},
to solve a free boundary problem of vortex dynamics that concerns steady-state two-dimensional  flow 
of a fluid in a wedge with impenetrable walls when the fluid contains a free boundary vortex domain, and there are no stagnation points on the walls.

A number of studies concern potential flows around  a vortex patch defined  \cite{Saffman1992} as a connected region of finite area containing uniform vorticity, surrounded by irrotational fluid.  
If the fluid inside a patch is at rest, then  the patch is  named as a hollow vortex.
Its boundary is a streamline, and pressure on the boundary is constant. 
In one of the first studies on this subject Pocklington \cite{{Pocklington1895}} solved
a two-dimensional problem of uniform flow of a fluid with a pair of hollow vortices of the same area and the circulations of opposite signs.
 The problem of a linear array of hollow vortices in an inviscid incompressible fluid when each vortex has a fore-and-aft symmetry was solved in   \cite{BakerSaffmanSheffield1976}. 
A model of flow around a hollow vortex placed in a corner was recently analyzed \cite{chr} by the method of 
conformal mappings based on the use of the Schottky-Klein prime function \cite{Crowdy2020}.
Since  the solution to the problem
of a point vortex in a corner  \cite{chr}   has a stagnation point on each side of the corner, 
the model  of a hollow vortex  also admits a stagnation point on each side of the  corner.
Such a model naturally results in the presence of a streamline that is orthogonal to the walls at the stagnation points and
separates the flow into two not mixing regions, a bounded corner region with  the hollow vortex inside
and the external region.

Sadovskii  \cite{Sadovskii1971}  considered the problem of splicing of a vortex flow symmetric  with respect to the 
$x$-axis in a potential flow
without rigid boundaries when the length (not the actual profile) of the vortex domain is prescribed, 
$-1\le x\le 1$, with the vortex distribution
$\Go(x,y)=-\GO_0\sgn y$, $\GO_0$ is a positive constant.  Two nonlinear equations for the vortex sheet strength
and the vortex profile function $y=f(x)$ were derived and solved numerically. The solution  \cite{Sadovskii1971}
shows that the streamline along the $x$-axis branches at the points $x=\pm 1$ and forms 
two cusps at the branch points.  
A Batchelor flow in a half-plane around a Sadovsky vortex
attached to the wall and 
a rotational right-angle corner flow 
was examined in
 \cite{MooreSaffmanTanveer1988}. In both cases the flow is separated by a vortex sheet 
 that is also a dividing streamline. The separated streamline is attached smoothly to the wall forming cusps at the attachment points. Using a Fourier series representation of the shape with a cuspidal
 behavior built in they recovered numerically the flow characteristics in both divided flow regions.

In our model of the flow in a wedge around a Sadovskii vortex (Fig 1a), it is assumed that the  velocity magnitude equals a constant,  $U$, on the whole vortex  boundary. In general, the tangential velocity is discontinuous across the vortex boundary, and the boundary itself is a vortex sheet.
Similarly to  \cite{Sadovskii1971}, in the vortex flow domain, there is 
a line (Fig. 1b) that separates the vortex domain into two subdomains not necessarily of the same area
with opposite circulations.
We do not determine the flow inside the vortex domain. On 
assuming that the constant speed   $U$ and the circulation, $\GG$, around the vortex boundary that is a branched streamline  are prescribed and nonzero we determine 
the vortex boundary.
The problem of the irrotational flow in the wedge exterior to the vortex is exactly solved
 by constructing the conformal mapping from the exterior of two cuts, $l_0=[m,\infty)$
 and $l_1=[0,1]$, onto the flow domain.
 Here, $m$ is a free parameter, and $m>1$.
The semi-infinite cut  is mapped onto the boundary of the wedge, the wedge vertex is the image of
some point $a$ on the lower side of the contour $l_0$, whilst the finite cut  is mapped onto the boundary of the vortex domain.
The conformal mapping is presented in terms of two functions. One of them is a rational function on an
elliptic surface $\CR$ topologically equivalent to a torus. The second one is the solution to a symmetric
 Riemann-Hilbert problem   on the same  Riemann surface $\CR$. 
Solutions to
 Riemann-Hilbert problems on genus-0, 1 and 2 Riemann surfaces were used 
  \cite{AntipovSilvestrov2007, AntipovSilvestrov2008, az, AntipovZemlyanova2009b, ZemlyanovaAntipov2012}
 to reconstruct the conformal mappings associated with
 supercavitating flows in simply, doubly, and triply-connected domains. It turns out that the Weierstrass analogue 
 of the Cauchy kernel \cite{zve}
 employed in these studies of supercavitating flows is not applicable to the Riemann-Hilbert problem arising
 in the vortex problem in a wedge. The reason  is that the Weierstrass  kernel is not sufficiently 
 decaying at infinity with respect to the variable of integration, and its use leads to divergent integrals.
 We propose a new analogue of the Cauchy kernel on a hyperelliptic surface with the properties needed. 
In the elliptic case,  this kernel is bounded at infinity and has simple poles at two bounded points of the surface.
These poles generate unacceptable essential singularities of the factorization function. They are removed by solving
the associated Jacobi inversion problem.
The factorization of the coefficient of the Riemann-Hilbert problem is  further used to derive the general solution to the Riemann-Hilbert problem that possesses a rational function on the surface $\CR$.
By satisfying  the additional conditions we arrive at a real transcendental equation
with respect to the preimage $a$ of the wedge vertex and the parameter $m$. It turn out that, when the infinite point 
of the cut $l_0$ falls into the infinite point of the wedge, the parameter $a=m$. This gives 
rise to a  family of conformal mappings possessing  the free geometric parameter $m$
and  two problem parameters, the wedge angle $\Ga$ and $\GG/U$.

This paper is organized as follows. In Section 2, we state the problem of a Sadovskii vortex in a wedge
and express the conformal map through two functions, $\Go_0(\Gz)$ and  $\Go_1(\Gz)$.
In Section 3, we determine the  function $\Go_0(\Gz)$. We derive  and solve the Riemann-Hilbert problem
on the elliptic surface $\CR$ for the function 
$\Go_1(\Gz)$ in Section 4. The exact formula for the conformal map is presented in Section 5.
We employ this map to reconstruct the vortex  and the wedge boundaries and report 
the numerical results in Section 6.

\setcounter{equation}{0}

\section{Model problem formulation and a conformal mapping}\label{s2.1}

\begin{figure}[t]
\centerline{
\scalebox{0.7}{\includegraphics{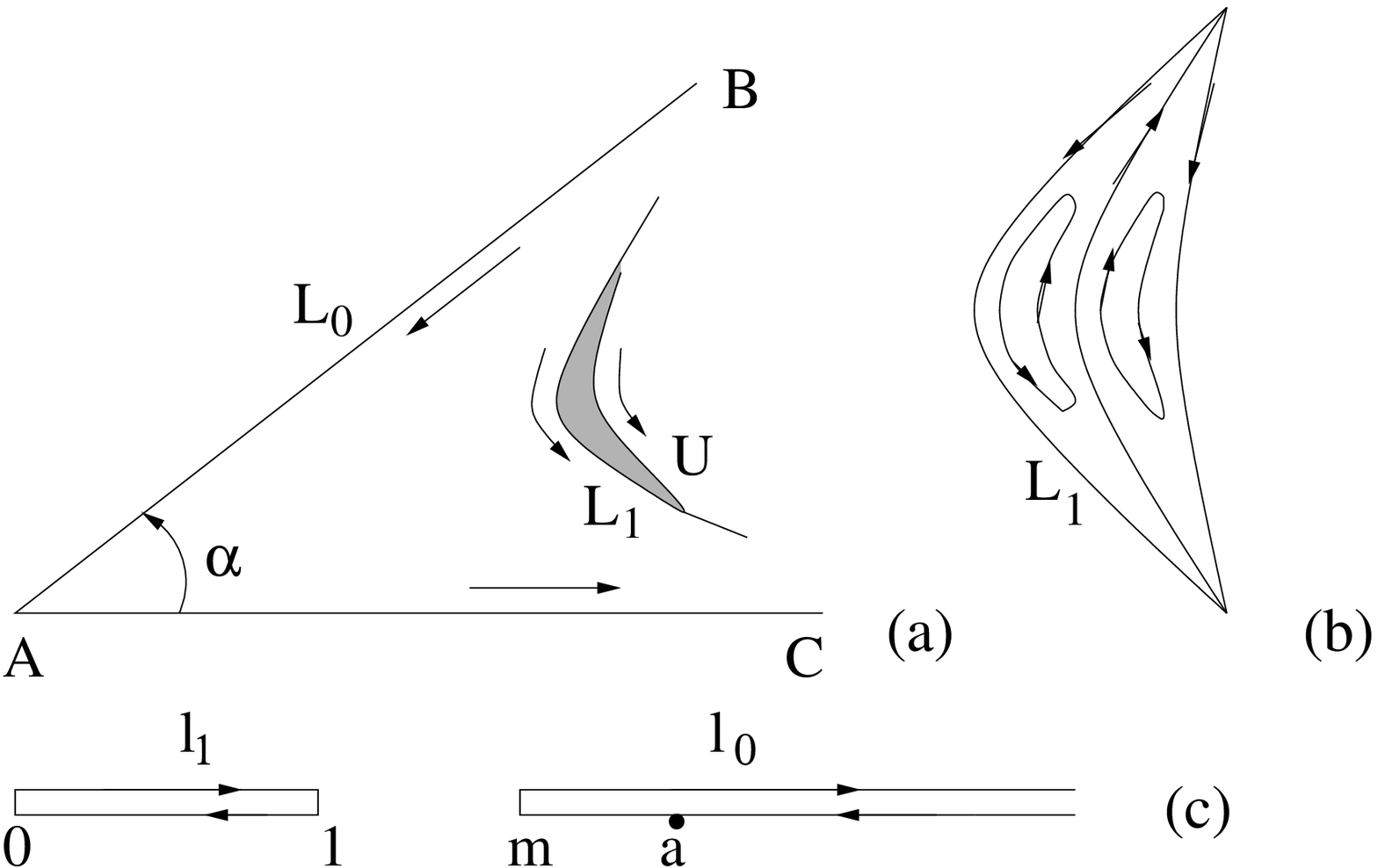}}}
\caption{(a): Potential flow in a wedge around a Sadovskii vortex. (b): Mechanism of flow in the vortex. (c):  Preimages of the vortex and the wedge boundaries in the parametric plane.}
\label{fig1}
\end{figure} 

We consider the following two-dimensional problem of fluid dynamics.

{\sl Let  an inviscid incompressible fluid flow in a wedge
$\CW=\{z\in{\Bbb C}: 0<|z|<\infty, 0<\arg z<\Ga\}$ whose  boundary $L_0$ is formed by rigid impenetrable walls 
$AC=\{z\in{\Bbb C}: z=r>0\}$ and $AB=\{z\in{\Bbb C}: z=re^{i\Ga}, r>0, 0<\Ga\le 2\pi\}$. 
Suppose that in the flow domain there is a vortex $\CP$  (Fig.1a) whose boundary $L_1$ 
and exact location  are unknown {\it a priori}.
Fluid speed is assumed to be constant  everywhere in the vortex boundary,
$|\Bv|=U$, $z\in L_1$, the circulation around the vortex is $\GG$, no stagnation points on the walls are admitted,
and the flow exterior to the vortex domain is irrotational. It is required to determine the vortex boundary and the flow
exterior to the vortex.}

The two-dimensional flow in the domain $\tilde\CD=\CW\setminus \CP$ is characterized by a complex potential $w(z)=\Gf(x,y)+i\psi(x,y)$
whose derivative describes the velocity vector $\Bv=(u_x,u_y)$
by $dw/dz=u_x-iu_y$, $z\in\tilde \CD$. Here, $\Gf$ is the velocity potential and
$\psi$ is the stream function. The potential $w(z)$ is an analytic multi-valued function in the flow
domain with a  cyclic period $\GG$,
\beq
\int_{L_1}\fr{dw}{dz}dz=\GG,
\label{2.0}
\eeq
and satisfies the following boundary conditions.  On the wedge walls they are
\beq
\arg\fr{dw}{dz}=\left\{
\begin{array}{cc}
\pi-\Ga, & z\in AB,\\
0, & z\in AC,\\
\end{array}
\right.
\quad \I w=c_0=\const, \quad z\in L_0.
\label{2.1}
\eeq
These conditions mean that our model excludes the presence of stagnation points on the upper and lower boundaries of the wedge.
Otherwise, if they were enforced, then $\arg dw/dz$ would be a piecewise constant function on the boundary.
For example, if $r=b\in(0,+\infty)$ were a stagnation point on the lower boundary  $AC$ and the velocity direction as $r\to\infty$
is the same  as in Fig.1,
then $\arg dw/dz$ would be $0$ for $r>b$ and 
$\pi$ for $0<r<b$. 

The boundary $L_1$ of the vortex is a streamline, $\psi(x,y)=c_1=\const$, 
$z\in L_1$, and since the flow speed $|dw/dz|$ is a constant $U$, 
 the boundary conditions on the vortex boundary read
\beq
\I w=c_1, \quad \left|\fr{dw}{dz}\right|=U, \quad z\in L_1.
\label{2.2}
\eeq
At large distances,  the complex potential behaves as
\beq
w(z)\sim \fr{\GG}{2\pi i}\log z+C z^\Gb, \quad z\to\infty, \quad z\in\tilde \CD,
\label{2.3}
\eeq
where $\Gb>0$ is a parameter to be determined,
and $C$ is a constant.

The vortex flow inside the vortex domain is 
characterized similarly to the Sadovskii model  \cite{Sadovskii1971}. A streamline
branches at some interior point in the wedge, the two branches merge at another point, and a finite region
is formed. At these two points,
the streamline branches form cusps. 
Otherwise, if the cusps are excluded,  two inadmissible stagnation points in the vortex boundary have to be enforced. The presence of the cusps will be also confirmed numerically. The vortex domain itself is divided into
two subdomains with vortex distributions of different sign (Fig. 1b).
 Apart from  the condition (\ref{2.0})
that prescribes the circulation around the vortex, owing to the second condition in (\ref{2.2})
and  the Bernoulli's equation, we know that the pressure distribution is constant on the vortex boundary.
In addition, we remark that the conditions  (\ref{2.2}) do not rule out a possible discontinuity
of the tangential velocity across the vortex boundary. In this case the boundary is a vortex sheet, and
as in \cite{Sadovskii1971}, the Bernoulli constant jumps across the vortex boundary.

To deal with the problem, consider another complex plane, a $\Gz$-plane, cut along two segments, 
$l_0=[m,+\infty)$ ($m>1$) and $l_1=[0,1]$ (Fig. 1c).
Denote the exterior of these two slits by $\CD$.  
 Let $z=f(\Gz)$ be a conformal map of this doubly connected domain onto the flow domain
 $\tilde\CD$. We assume that the function $f(\Gz)$ maps the slits $l_0$ and $l_1$ onto
 the wedge and vortex boundaries $L_0$ and $L_1$, respectively. The positive directions of the contours $l_0$ and $l_1$ are chosen to be consistent with those in the physical domain — when a point $\Gz$ traverses the contours,
 the domain $\CD$ is seen on the left. Since the direction of $l_1$ corresponds to the negative
 direction of the contour $L_1$ with respect to the interior of the vortex domain, the circulation
  $\GG<0$.
 For the upper and lower sides of the contours 
 $l_0$ and $l_1$, we shall use the notations $l_0^\pm$ and $l_1^\pm$, respectively.
 We choose the map such that the infinite point $\Gz=\infty\in\CD$ falls into the infinite point of the flow
 domain,
 $z=\infty\in \tilde\CD$. Without loss, it is assumed that the preimage of the wedge vertex
 $A$ is a point $a$ lying in the lower side of the semi-infinite cut including the point $m$.
 This assumption inevitably gives us the preimages of the points $B$ and $C$
 in the parametric plane: they are $+\infty-i0$ and $+\infty+i0$, respectively.
 
In order to reconstruct the conformal mapping, we express the derivative
$f'(\Gz)$ through the derivative $dw/d\Gz$ and the logarithmic hodograph 
variable as
\beq
\fr{df}{d\Gz}=\Go_0(\Gz)e^{-\Go_1(\Gz)},
\label{2.4}
\eeq
where
\beq
\Go_0(\Gz)=\fr{dw}{Ud\Gz}, \quad
\Go_1(\Gz)=\log\fr{dw}{Udz}.
\label{2.5}
\eeq

\setcounter{equation}{0}

\section{Function $\Go_0(\Gz)$ and asymptotics of the map}\label{S3}

Since the imaginary part of the function $w(z)$ is 
constant on $L_0$ and $L_1$, we have $\I dw/ d\Gz=0$ on the contours $l_0$ and $l_1$.
In virtue of (\ref{2.5})  the function $\Go_0(\Gz)$ is analytic in the domain $\CD={\Bbb C}\setminus(l_0\cup l_1)$, satisfies the boundary condition
\beq
\I \Go_0(\Gx)=0, \quad \Gx\in l_0\cup l_1,
\label{3.1}
\eeq
and may have integrable singularities at the finite endpoints of the contours  $l_0$ and $l_1$.
At infinity, $\Go_0(\Gz)=O(\Gz^{-\Gm})$ with $\Gm>0$.
The most general form of such a function is
\beq
\Go_0(\Gz)=\fr{i(N_0+N_1\Gz)}{p^{1/2}(\Gz)},
\label{3.2}
\eeq
where $N_0$ and $N_1$ are arbitrary real constants, and
\beq
p(\Gz)=\Gz(1-\Gz)(\Gz-m).
\label{3.3}
\eeq
The function $p^{1/2}(\Gz)$ is analytic in the $\Gz$-plane cut along the lines $l_0$
and $l_1$.
Its single branch is
fixed by the condition $p^{1/2}(\Gz)=i\sqrt{|p(\Gx)|}$
as $\Gz=\Gx+i0$, $\Gx>m$. This branch has the following properties:
$$
p^{1/2}(\Gz)=\pm i\sqrt{|p(\Gx)|}, \quad \Gz=\Gx\pm i0, \quad m<\Gx<+\infty,
$$
$$
p^{1/2}(\Gz)=-\sqrt{|p(\Gx)|}, \quad \Gz=\Gx, \quad 1<\Gx<m,
$$
$$
p^{1/2}(\Gz)=\mp i\sqrt{|p(\Gx)|}, \quad \Gz=\Gx\pm i0, \quad 0<\Gx<1,
$$
\beq
p^{1/2}(\Gz)=\sqrt{|p(\Gx)|}, \quad \Gz=\Gx, \quad -\infty<\Gx<0.
\label{3.4}
\eeq
The function $\Go_0(\Gz)$ is bounded as $\Gz\to a$ if $a\ne m$
and has the square root singularity $\Go_0(\Gz)=O((\Gz-a)^{-1/2})$ otherwise.

We next analyze the behavior of the function $f(\Gz)$  at the 
points $\Gz=a$ and $\Gz=\infty.$
Consider a copy of the $\Gz$-plane cut along the semi-infinite
segment $l_0=[m,\infty)$. By the map $f_1(\Gz)=\sqrt{\Gz-m}$ it is transformed into
the upper half-plane with the points $m$ and $a\in l_0^-$ being transformed
into the origin and the point $-\sqrt{a-m}$, respectively. By the translation transformation
$f_2(\Gz)=f_1(\Gz)+\sqrt{a-m}$ we move the point $-\sqrt{a-m}$ into the origin. Finally, by implementing the transformation
$z=f_2^{\Ga/\pi}(\Gz)$, we map the upper half-plane into the wedge $W$ in the 
$z$-plane. The composition of these three maps has the form
\beq
z=(\sqrt{\Gz-m}+\sqrt{a-m})^{\Ga/\pi}=\fr{(\Gz-a)^{\Ga/\pi}}{(\sqrt{\Gz-m}-\sqrt{a-m})^{\Ga/\pi}}.
\label{3.5}
\eeq
Note that this map transforms the finite cut $l_1$ into a closed contour inside the wedge $W$. The map (\ref{3.5}) and the function $f(\Gz)$ share the same asymptotics
at the points $a$ and $\infty$.
Consequently, if $a\ne m$, then formula (\ref{3.5}) implies
$$
f(\Gz)\sim K_0(\Gz-a)^{\Ga/\pi}, \quad \Gz\to a\in l_0^-,
\quad
f(\Gz)\sim K_1, \quad \Gz\to \bar a\in l_0^+,
$$
\beq
f(\Gz)\sim \Gz^{\Ga/(2\pi)}, \quad \Gz\to \infty.
\label{3.6}
\eeq
Here, $K_0$ and $K_1$ are some constants. 

In the case $a=m$, the asymptotics of the function $f(\Gz) $ at infinity is the same, whilst
at the left endpoint of the cut $l_0$, we have
 \beq
f(\Gz)\sim K_0(\Gz-a)^{\Ga/(2\pi)}, \quad \Gz\to a=m.
\label{3.7}
\eeq

\setcounter{equation}{0}

\section{Function $\Go_1(\Gz)$ and a Riemann-Hilbert problem on a torus}

\subsection{Hilbert problem for the function $\Go_1(\Gz)$}

We now turn our attention to the function $\Go_1(\Gz)$. Comparing
the asymptotics (\ref{2.3}) and (\ref{3.6}) at infinity and using formula (\ref{3.2}) we deduce from
(\ref{2.4}) and (\ref{2.5}) 
\beq
\Gb=\fr{\pi}{\Ga}, \quad \fr{dw}{Udz}\sim\fr{\pi C}{\Ga}\Gz^{(\pi-\Ga)/(2\pi)}, \quad \Gz\to\infty,
\label{4.1}
\eeq
and therefore
\beq
\Go_1(\Gz)\sim \fr{\pi-\Ga}{2\pi}\log\Gz, \quad \Gz\to\infty.
\label{4.2}
\eeq

In a similar fashion we determine the asymptotics of the function $\Go_1(\Gz)$ at the point $a$.
If $a\ne m$ ($a\in l_0^-$), then we have 
\beq
\fr{dw}{Udz}\sim C'(\Gz-a)^{1-\Ga/\pi}, \quad \Gz\to a, \quad 
\fr{dw}{Udz}\sim C'', \quad \Gz\to \bar a,
\label{4.3}
\eeq
where  $C'$ and $C''$ are constants. The second formula in (\ref{2.5}) yields
\beq
\Go_1(\Gz)\sim\fr{\pi-\Ga}{\pi}\log(\Gz-a), \quad \Gz\to a, \quad
\Go_1(\Gz)\sim \log C'', \quad \Gz\to \bar a.
\label{4.4}
\eeq
If it happens that $a=m$, then
\beq
\Go_1(\Gz)\sim\fr{\pi-\Ga}{2\pi}\log(\Gz-a), \quad \Gz\to a.
\label{4.5}
\eeq
At the point $\bar a$ the function $\Go_1(\Gz)$ is bounded and at infinity has the same asymptotics
(\ref{4.2}) as in the case $a\ne m$.

Analyze next the boundary conditions. Observe in equations (\ref{2.1}) and (\ref{2.2}) that the values of
the argument and modulus of the function $dw/dz$  give rise to the imaginary and real parts 
 of  the function $\Go_1(\Gz)$. We have
 \beq
 \I\Go_1(\Gz)=\left\{
 \begin{array}{cc}
 \pi-\Ga, &  \Gz\in l_0',\\
 0,  & \Gz\in l_0'',\\
 \end{array}
 \right.\quad
 \R\Go_1(\Gz)=0, \quad \Gz\in l_1.
 \label{4.6}
 \eeq
 Here, $l_0'=[a,+\infty)^-$,  $l_0''=[m,a]^-\cup[m,+\infty)^+$, and the superscripts
 $+$ and $-$ indicate that the segments belong to the upper and lower
 sides of the contour $l_0$, correspondingly.
 The problem for the function $\Go_1(\Gz)$ obtained is classified as a Hilbert
 problem of the theory of analytic functions on two two-sided positively oriented contours
 (the exterior of the loops $l_0$ and $l_1$ is on the left when a point $\Gz$ traverses the contours
 in the positive direction). 
 To solve this problem, we reduce it to a Riemann-Hilbert problem on two contours
 of a symmetric genus-1 Riemann surface  topologically equivalent to a torus.

 \subsection{Riemann-Hilbert problem for an auxiliary function $\GF(\Gz,u)$ on a torus}
 
 Let  ${\Bbb C}_1$ and  ${\Bbb C}_2$ be two copies of the extended complex $\Gz$-plane 
cut along the segments $l_1=[0,1]$ and $l_0=[m,+\infty)$ and
$\CR$ be a genus-1 Riemann surface defined by the algebraic equation
\beq
u^2=p(\Gz),\quad  p(\Gz)=\Gz(1-\Gz)(\Gz-m).
\label{4.7}
\eeq
The surface is formed by gluing the two sheets  together in such a way that
the upper sides $l_j^+$
of the cuts $l_j\subset {\Bbb C}_1$ are glued to the lower 
sides $l_j^-$ of the cuts $l_j\subset {\Bbb C}_2$, and the sides
$l_j^-\subset {\Bbb C}_1$ are glued to $l_j^+\subset {\Bbb C}_2$ ($j=0,1$).
Let  $p^{1/2}(\Gz)$ be the branch fixed in Section \ref{S3}.
Then the function $u(\Gz)$  
\beq
u=\left\{\begin{array}{cc}
p^{1/2}(\Gz),  & (\Gz,u)\in{\Bbb C}_1,\\
-p^{1/2}(\Gz),  & (\Gz,u)\in{\Bbb C}_2,\\
\end{array}
\right.
\label{4.8}
\eeq
is uniquely defined and single-valued on $\CR$.
The pairs $(\Gz, p^{1/2}(\Gz))$ and $(\Gz, -p^{1/2}(\Gz))$ correspond to the points with affix $\Gz$ lying
on the upper and lower sheets, ${\Bbb C}_1$ and ${\Bbb C}_2$, respectively. 
The contour $\CL=l_0\cup l_1$
is the symmetry line for the elliptic surface $\CR$
which splits the surface into two symmetric halves with symmetric points 
$(\Gz,u)\in{\Bbb C}_1$ and $(\Gz_*,u_*)=(\bar\Gz, -u(\bar\Gz))\in{\Bbb C}_2$.

Introduce an auxiliary  function on the surface  $\CR$
\beq
\GF(\Gz,u)=\left\{
\begin{array}{cc}
-i\Go_1(\Gz), & (\Gz,u)\in{\Bbb C}_1,\\
i\,\overline{\Go_1(\bar\Gz)}, & (\Gz,u)\in{\Bbb C}_2.\\
\end{array}
\right.
\label{4.9}
\eeq
Everywhere on the torus this function satisfies the symmetry condition
\beq
\overline{\GF(\Gz_*,u_*)}=\GF(\Gz,u),\quad(\Gz,u)\in\CR.
\label{4.10}
\eeq
On the symmetry line $\CL$, its boundary values are expressed through the real and imaginary
values of the function $\Go_1(\Gx)$,
$$
\GF^+(\Gx,v)=-i\Go_1(\Gx)=-i\R\Go_1(\Gx)+\I\Go_1(\Gx),
$$
\beq
\GF^-(\Gx,v)=i\overline{\Go_1(\Gx)}=i\R\Go_1(\Gx)+\I\Go_1(\Gx),\quad 
(\Gx,v)\in \CL.
\label{4.11}
\eeq
Here and henceforth, $v=u(\Gx)$, $\GF^+(\Gx,v)$  and $\GF^-(\Gx,v)$
are the limit values of the function $\GF(\Gz,u)$ on the contour
$\CL$ from the upper  and lower sheet of the surface, respectively.
On summarizing the properties of the function $\GF(\Gz,u)$ we state the following Riemann-Hilbert problem.

{\sl Find all symmetric functions $\GF(\Gz,u)$, $\overline{\GF(\Gz_*,u_*)}=\GF(\Gz,u)$, 
$(\Gz,u)\in\CR$,
 analytic in $\CR\setminus \CL$,
H\"older-continuous up to the boundary $\CL$ apart from the singular points
$\Gz=a$ and $\Gz=\infty$
with the boundary values satisfying the relation
\beq
\GF^+(\Gx,v)=G(\Gx,v)\GF^-(\Gx,v)+g(\Gx,v), \quad (\Gx,v)\in\CL,
\label{4.12}
\eeq
 Here, 
 \beq
G(\Gx,v)=\left\{
\begin{array}{cc}
-1, & (\Gx,v)\in l_0\\
1, & (\Gx,v)\in l_1,\\
\end{array}
\right.
\quad
g(\Gx,v)=\left\{
\begin{array}{cc}
2(\pi-\Ga), & (\Gx,v)\in l_0',\\
0, & (\Gx,v)\in l_0'',\\
0, & (\Gx,v)\in l_1.\\
\end{array}
\right.
\label{4.13}
\eeq
The function $\GF(\Gz,u)$ has a logarithmic singularity at the
infinite point  $(\infty,\infty)$,
\beq
\GF(\Gz,u)\sim\fr{\pi-\Ga}{2\pi i}\log\Gz, \quad \Gz\to \infty, 
\label{4.14}
\eeq
It is bounded at the point $(\bar a,u(\bar a))$ and has a logarithmic singularity at the point
$(a, u(a))$,
\beq
\GF(\Gz,u)\sim\fr{\pi-\Ga}{\Gs\pi i}\log(\Gz-a),\quad \Gz\to a, 
\label{4.15}
\eeq
where $\Gs=1$ if $a\ne m$ and $\Gs=2$ otherwise.}

 \subsection{Analogue of the Cauchy kernel}
 
The standard solution device for scalar Riemann-Hilbert problems
on elliptic and hyperelliptic  surfaces is the Wierstrass kernel \cite{zve}, an analogue  
of the Cauchy kernel for Riemann surfaces. It has the form
\beq
dW=\fr12\left(1+\fr{u}{v}\right)\fr{d\Gx}{\Gx-\Gz}, \quad u=u(\Gz), \quad v=u(\Gx). 
\label{4.16}
\eeq
This kernel is applicable if the Riemann-Hilbert problem contour is bounded
or the density decays at infinity at a sufficient rate otherwise.
When the contour  is unbounded
and the density does not decay at infinity fast enough, the employment of the kernel (\ref{4.16}) gives
rise to divergent integrals.
In our case the contour $\CL$ of the Riemann-Hilbert problem comprises a finite contour $l_1$ and a semi-infinite contour $l_0$.
Since the function $g(\Gx,v)$ defined in (\ref{4.13}) is a nonzero constant in the contour
$l_0'=[a,\infty)^-$, we need a kernel different from the one in (\ref{4.16}).
 
We propose to use a new analogue of the Cauchy kernel convenient for our purposes. 
Consider the differential
\beq
dV=\fr12\left[1+\fr{u}{v}\fr{(\Gx-\Gz_0)^2}{(\Gz-\Gz_0)^2}\right]
\left(\fr{1}{\Gx-\Gz}-\fr{1}{\Gx-\Gz_0}\right)d\Gx,
\label{4.17}
\eeq
where $\Gz_0$ is an arbitrary real fixed point not lying on the contours $l_0$ and
$l_1$ (in our case it is convenient to select $\Gz_0$ by the condition $\Gz_0<0$). 
Equivalently, the kernel $dV$ may be written in the form
\beq
dV=\fr12\left[\fr{\Gz-\Gz_0}{\Gx-\Gz_0}+\fr{u}{v}\fr{\Gx-\Gz_0}{\Gz-\Gz_0}\right]
\fr{d\Gx}{\Gx-\Gz}.
\label{4.17'}
\eeq
This kernel possesses the following properties.

(i) If $\Gz$ is a fixed bounded point on the torus $\CR$ and $\Gx\to\pm i0+\infty$, then
$dV$ is decaying as $\Gx^{-3/2}$ that is
\beq
dV\sim \pm\fr{1}{2i}\fr{u(\Gz)}{\Gz-\Gz_0}\Gx^{-3/2}d\Gx, \quad \Gx\to\pm i0+\infty.
\label{4.18}
\eeq
 
 (ii) If $\Gx$ is a fixed point lying either in the contour $l_1$ or any finite part of the contour $l_0$
 and $\Gz\to\infty$, 
 then the kernel $dV$ is bounded,
\beq
dV\sim \left[-\fr{1}{2(\Gx-\Gz_0)}+O(\Gz^{-1/2})\right]d\Gx, \quad \Gz\to \infty.
\label{4.19}
\eeq
  
(iii) If $\Gz\to\Gx\in\CL$, then the kernel behaves as the Cauchy kernel.
Indeed, if $\Gz\in{\Bbb C}_1$ and $\Gz\to\Gx$, then for $\Gx\in \CL\subset{\Bbb C}_1$
\beq
dV\sim \fr{d\Gx}{\Gx-\Gz}, 
\label{4.20}
\eeq
whilst for  $\Gx\in \CL\subset{\Bbb C}_2$, $u/v\sim -1$, and the kernel $dV$
is bounded. The analysis is similar when  $\Gz\in{\Bbb C}_2$.
 
(iv) Since $\Gz_0$ is a real number, the kernel $dV$ is symmetric  with respect to the symmetry contour $\CL$,
\beq
dV((\Gx,v),(\Gz,u(\Gz)))=\ov{dV((\Gx,v),(\bar\Gz,-u(\bar \Gz)))}.
\label{4.21}
\eeq 
 
 (v) At the two points of the torus with affix $\Gz_0$, $(\Gz_0, u(\Gz_0))$ and
 $(\Gz_0, -u(\Gz_0))$, the kernel  $dV$ has simple poles. They will give rise to
 essential singularities of the factorization function associated with the Riemann-Hilbert problem (\ref{4.12})
and will have to be removed by solving a Jacobi inversion problem.

The kernel we presented can be generalized in the case of a hyperelliptic surface of any finite genus $\Gr$.
We write down such an analogue since similar problems of higher genera might arise in other applications. Let $\CR$ be the genus-$\Gr$ hyperelliptic surface
of the algebraic function $u^2=(\Gz-a_1)(\Gz-a_2)\ldots (\Gz-a_{2\Gr+1})$ (the infinite point is a branch point)
or $u^2=(\Gz-a_1)(\Gz-a_2)\ldots (\Gz-a_{2\Gr+2})$ (all the branch points are finite).
Then the counterpart of 
the kernel (\ref{4.17}) has the form
 \beq
dV=\fr12\left[1+\fr{u}{v}\prod_{j=0}^{\Gr}\fr{\Gx-\Gz_j}{\Gz-\Gz_j}\right]
\left(\fr{1}{\Gx-\Gz}-\fr{1}{\Gx-\Gz_0}\right)d\Gx.
\label{4.21.0}
\eeq 
Here, $\Gz_0,\ldots,\Gz_{\Gr}$ are arbitrary fixed not necessarily distinct points in the real axis not falling on the 
 contour of the Riemann-Hilbert problem that is the union of the branch cuts of the Riemann surface.

 \subsection{Factorization problem}

 \begin{figure}[t]
\centerline{
\scalebox{0.7}{\includegraphics{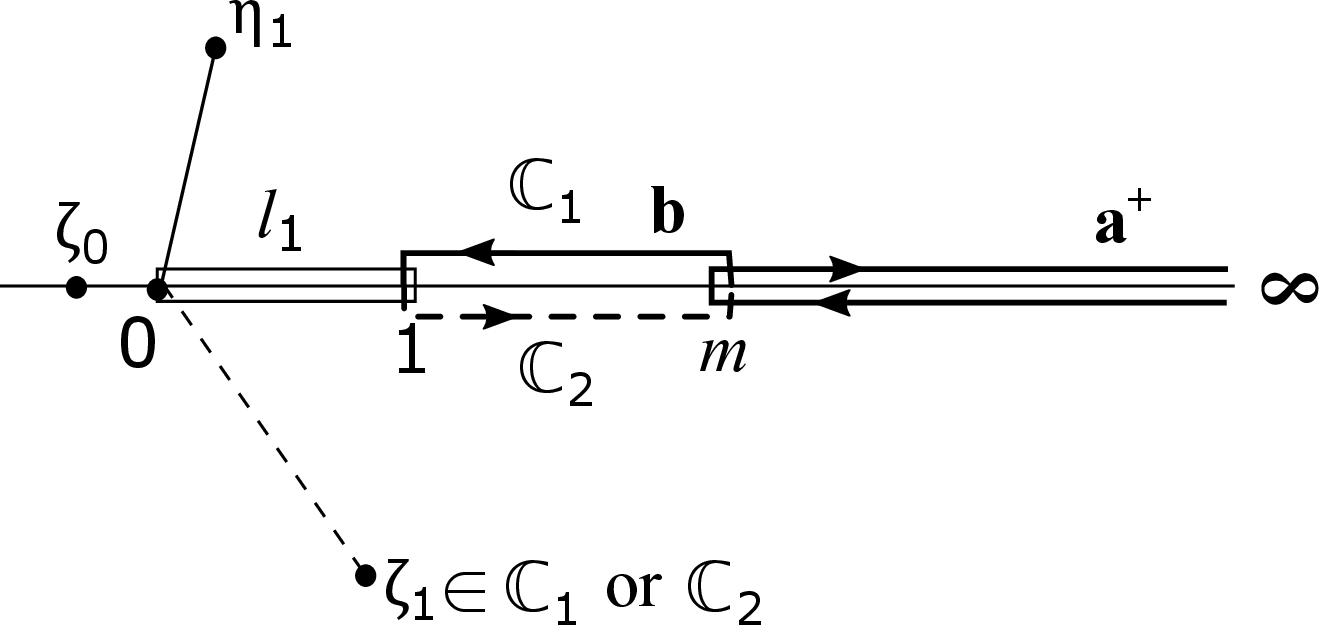}}}
\caption{Canonical cross-sections and the contour $\Gg=\Bp_1\Bq_1$,  
$\Bp_1=(\Gn_1,u(\Gn_1))\in {\Bbb C_1}$, and  $\Bq_1=(\Gz_1,u(\Gz_1))\in\CR$.}
\label{fig2}
\end{figure} 

We aim to determine a piece-wise meromorphic function $X(\Gz,u)$ symmetric on the torus,
\beq
X(\Gz,u)=\ov{X(\bar\Gz,-u(\bar\Gz))}, \quad (\Gz,u)\in \CR\setminus \CL,
\label{4.22.0}
\eeq 
continuous through the contour $l_1$,
discontinuous through the contour $l_0$,
and whose boundary values on the contour satisfy the equation
\beq
X^+(\Gx,v)=G(\Gx,v){X^-(\Gx,v)},\quad (\Gx,v)\in \CL\subset\CR.
\label{4.22}
\eeq
Analyze the following function:
\beq
X(\Gz,u)=\Gc(\Gz,u)\Gc_1(\Gz,u)\ov{\Gc_1(\Gz_*,u_*)},
\label{4.23}
\eeq
where
$$
\Gc(\Gz,u)=\exp\left\{\fr{1}{2\pi i}\oint_{l_0}\log(-1) dV\right\},
$$
\beq
\Gc_1(\Gz,u)=\exp\left\{-\left(\int_{\Gg}
+n_a\oint\limits_{\rm\bf a}
+n_b\oint\limits_{\rm\bf b}\right)dV
\right\}.
\label{4.24}
\eeq
The contours $\Ba$ and  $\Bb$ constitute a system  of canonical cross-sections 
of the surface $\CR$. They are chosen as follows. The contour  {\bf a}
lies on both sheets
of the surface and coincides with the banks of the semi-infinite cut $l_0$  (Fig.2). 
The loop {\bf b} consists of the segments $[m,1]\subset{\Bbb C}_1$
and $[1,m]\subset{\Bbb C}_2$.  
The positive direction on the loop {\bf a} is chosen such that the upper sheet is on the left.
The loop {\bf a} intersects the loop {\bf b} at the branch point $\Gz=m$ from left to the right.

The contour $\Gg$ is a continuous curve whose starting and terminal points are $\Bp_1=(\eta_1, u(\eta_1))$
and $\Bq_1=(\Gz_1,u(\Gz_1))$, respectively. 
The point $\Bp_1$ is an arbitrary fixed point lying in either sheet of the surface.
Without loss, select  $\Bp_1=(\Gn_1,p^{1/2}(\Gn_1))\in{\Bbb C_1}$.
As for the point $\Bq_1$,  it has to be determined and may fall 
on either sheet. The integers $n_a$ and $n_b$ cannot be fixed {\it a priori} either
and are to be recovered from a condition which guarantees the boundedness   
of the solution to the problem (\ref{4.22}) at the two points of the torus with affix $\Gz_0$, 
the simple poles of the kernel $dV$.
The contour $\Gg$  is chosen such that it does not cross the {\bf a}- and {\bf b}-cross-sections. In the case $\Bq_1\in{\Bbb C}_2$,
it passes through the point $\Gz=0$, a branch point of the surface $\CR$, and consists of two parts,
$(\Bp_1,0)\subset {\Bbb C}_1$ and $(0,\Bq_1)\subset {\Bbb C}_2$.
If it turns out that the point $\Bq_1$ lies on the upper sheet, then the contour $\Gg$
can be taken as the straight segment joining these points provided it does not cross
the loops $l_0$ and $l_1$.  Otherwise, if the imaginary parts of the points $\Gn_1$ 
and $\Gz_1$ have different signs, both points $\Bp_1$ and $\Bq_1$ lie in the sheet
${\Bbb C_1}$, and the straight segment $\Bp_1 \Bq_1$  intersects the contour $\CL$, 
then $\Gg$ can be selected as a polygonal line consisting of two segments
$\Bp_1\Bp_0$ and $\Bp_0\Bq_1$ with 
 $\Bp_0=(\Gd, p^{1/2}(\Gd))$ and $\Gz_0<\Gd<0$.

Note that the integrals in (\ref{4.24}) have jumps through the contours $\Gg$, {\bf a}, and {\bf b},
respectively. However, since the jumps are multiples of $2\pi i$, the functions
in the right hand-side in (\ref{4.23}) are continuous through these contours. 
It is directly verified that the function (\ref{4.23}) satisfies the symmetry condition (\ref{4.22.0}).

The integral representation of the factorization function $X(\Gz,u)$  can be simplified. 
By employing the Cauchy theorem
we deduce  
\beq
\Gc(\Gz,u)=\exp
\left\{\fr{u(\Gz)}{4(\Gz-\Gz_0)}
\oint_{l_0} \fr{(\Gx-\Gz_0)d\Gx}{u(\Gx)(\Gx-\Gz)}\right\}.
\label{4.25}
\eeq
In the product $\Gc_1(\Gz,u)\ov{\Gc_1(\Gz_*,u_*)}$, the integrals over the cross-section $\Bb$
are cancelled and we have
$$
\Gc_1(\Gz,u)\ov{\Gc_1(\Gz_*,u_*)}=
\exp
\left\{-\fr12\int_\Gg\left(\fr{\Gz-\Gz_0}{\Gx-\Gz_0}+\fr{u(\Gz)}{u(\Gx)}\fr{\Gx-\Gz_0}{\Gz-\Gz_0}\right)
\fr{d\Gx}{\Gx-\Gz}
\right.
$$
\beq
\left.
-\fr12\int_\Gg\left(\fr{\Gz-\Gz_0}{\bar\Gx-\Gz_0}-\fr{u(\Gz)}{u(\bar\Gx)}\fr{\bar\Gx-\Gz_0}{\Gz-\Gz_0}\right)\fr{d\bar\Gx}{\bar\Gx-\Gz}-
\fr{n_a u(\Gz)}{\Gz-\Gz_0}
\oint_{\rm\bf a} \fr{(\Gx-\Gz_0)d\Gx}{u(\Gx)(\Gx-\Gz)}\right\}.
\label{4.25'}
\eeq
Since the contour $l_0$ and the loop $\Ba$ have the same direction and coincide pointwise, we  derive
the following representation for the solution of the factorization problem:
$$
X(\Gz,u)=\exp
\left\{\left(\fr12-2n_a\right)\fr{u(\Gz)}{i(\Gz-\Gz_0)}
\int_m^\infty \fr{(\Gx-\Gz_0)d\Gx}{\sqrt{|p(\Gx)|}(\Gx-\Gz)}
\right.
$$
\beq
\left.
-\fr12\int_\Gg\left(\fr{\Gz-\Gz_0}{\Gx-\Gz_0}+\fr{u(\Gz)}{u(\Gx)}\fr{\Gx-\Gz_0}{\Gz-\Gz_0}\right)
\fr{d\Gx}{\Gx-\Gz}
-\fr12\int_\Gg\left(\fr{\Gz-\Gz_0}{\bar\Gx-\Gz_0}-\fr{u(\Gz)}{u(\bar\Gx)}\fr{\bar\Gx-\Gz_0}{\Gz-\Gz_0}\right)\fr{d\bar\Gx}{\bar\Gx-\Gz}
\right\}.
\label{4.26}
\eeq
We describe the main properties of the function $X(\Gz,u)$ given by (\ref{4.26}).

(i) The function $X(\Gz,u)$ satisfies the symmetry condition (\ref{4.22.0}) and the boundary condition
(\ref{4.22}).

(ii) The limit values $X^+(\Gx\pm i0, u^\pm(\Gx))$ and $X^-(\Gx\pm i0, u^\pm(\Gx))$ on both sides of the contour $l_0$ are
expressed through the Cauchy principal value of the first integral in (\ref{4.26}) as
$$
X^+(\Gx\pm i0, u^\pm(\Gx))=i X(\Gx, u^\pm(\Gx)), \quad \Gx\in l_0^\pm.
$$
\beq
X^-(\Gx\pm i0, u^\pm(\Gx))=-i X(\Gx, u^\pm(\Gx)), \quad \Gx\in l_0^\pm.
\label{4.27}
\eeq

(iii) In virtue of the behavior of the Cauchy integral at the endpoints of the contour 
$\Gg$, the function $X(\Gz,u)$ has a simple zero at the point $\Bp_1\in{\Bbb C}_1$ 
and a simple pole at the point $\Bq_1\in\CR$.

(iv) As a consequence of the  poles  of the kernel at the points 
$(\Gz_0,u(\Gz_0))$ and $(\Gz_0,-u(\Gz_0))$, the function $X(\Gz,u)$ has essential 
singularities at these two points. 
They will be removed in the next section by solving the associated
Jacobi inversion problem.

\subsection{Jacobi inversion problem}

 To remove the essential singularities of the factorization function $X(\Gz,u)$
 at the points $(\Gz_0,\pm u(\Gz_0))$, we require
 \beq
  \mathop{\rm res}\limits_{\Gz=\Gz_0}\fr{u(\Gz)}{\Gz-\Gz_0}\left\{
 \left(\fr12-2n_a\right)
\int_m^\infty \fr{(\Gx-\Gz_0)d\Gx}{i\sqrt{|p(\Gx)|}(\Gx-\Gz)}
-\fr12\int_\Gg\fr{(\Gx-\Gz_0)d\Gx}{u(\Gx)(\Gx-\Gz)}
+\fr12\int_\Gg\fr{(\bar\Gx-\Gz_0)d\bar\Gx}{u(\bar\Gx)(\bar\Gx-\Gz)} \right\}=0.
\label{4.28}
\eeq
This brings us to the following nonlinear equation with respect to the point $\Bq_1$, the terminal point of the contour $\Gg$, and the integer $n_a$:
\beq
\left(\fr{1}{4}-n_a\right)
\oint\limits_{{\rm\bf a}}\fr{d\Gx}{u(\Gx)}-
\fr12\int\limits_{\Gg}\fr{d\Gx}{u(\Gx)}+\fr12\int\limits_{\Gg}\fr{d\bar\Gx}{\ov{u(\Gx)}}=0,
\label{4.29}
\eeq 
which can be considered \cite{az} as an imaginary analogue of the Jacobi inversion problem 
\beq
\int\limits_{0}^{\Bq_1}\fr{d\Gx}{u(\Gx)}+n_a\CA+n_b\CB=g_0,
\label{4.30}
\eeq
where  $\CA$ and $\CB$ are the $\CA$- and $\CB$- periods   of the abelian integral 
\beq
\int\limits_{0}^{\Bq}\fr{d\Gx}{u(\Gx)}
\label{4.31}
\eeq
given in terms of the Legendre complete elliptic integrals \cite{low}
$$
\CA=\oint\limits_{{\rm\bf a}}\fr{d\Gx}{u(\Gx)}=-2i\int\limits_m^\infty\fr{d\Gx}{\sqrt{|\Gx(\Gx-1)(\Gx-m)|}}=-4ik\BK(k),
$$
\beq
\CB=\oint\limits_{{\rm\bf b}}\fr{d\Gx}{u(\Gx)}=2\int\limits_1^m\fr{d\Gx}{\sqrt{|\Gx(\Gx-1)(\Gx-m)|}}=4k\BK'(k),
\label{4.32}
\eeq
$k=m^{-1/2}$, and
\beq
g_0=-ik\BK(k)+\int_0^{\Gn_1}\fr{d\Gx}{p^{1/2}(\Gx)}.
\label{4.33}
\eeq
The only one difference between the Jacobi problem (\ref{4.30}) and the one solved in \cite{az}
is the right hand-side $g_0$. By adjusting the formulas \cite{az} to our case
we can describe the solution procedure as follows.
Compute the affix of the point $\Bq_1$ by
\beq
\Gz_1={\rm sn}^2\fr{i g_0}{2k}
\label{4.34}
\eeq
and evaluate the integrals
\beq
I_{\pm}=\int\limits_{0}^{\Gn_1}
\fr{d\Gx}{p^{1/2}(\Gx)}\pm\int\limits_{0}^{\zeta_1}\fr{d\Gx}{p^{1/2}(\Gx)}-ik\BK(k)
\label{4.35}
\eeq
and the numbers $n_a$ and $n_b$
\beq
n_a=-\fr{\I I_-}{4k\BK(k)}, \quad n_b=\fr{\R I_-}{4k\BK'(k)}.
\label{4.36}
\eeq
If both of the numbers $n_a$ and $n_b$ are integers, then 
the point $\Bq_1$ lies in the upper sheet ${\Bbb C}_1$, $ \Bq_1=(\Gz_1, p^{1/2}(\Gz_1))$.
Otherwise, $\Bq_1\in {\Bbb C}_2$, $ \Bq_1=(\Gz_1, -p^{1/2}(\Gz_1))$, and
 the numbers $n_a$ and $n_b$ computed by the formulas 
 \beq
n_a=-\fr{\I I_+}{4k\BK(k)}, \quad n_b=\fr{\R I_+}{4k\BK'(k)}
\label{4.37}
\eeq
are integers.
If the point $\Bq_1$ and the integer $n_a$ are determined according to this procedure,
then the function $X(\Gz,u)$ is bounded at the points $(\Gz_0,\pm u(\Gz_0))$ and provides
a solution to the factorization problem (\ref{4.22.0}), (\ref{4.22}).

\subsection{Solution to the Riemann-Hilbert problem}

First we use the factorization of the function $G(\Gx,v)$ and 
rewrite the boundary condition  (\ref{4.12}) in the form
\beq
\fr{\GF^+(\Gx,v)}{X^+(\Gx,v)}=\fr{\GF^-(\Gx,v)}{X^-(\Gx,v)}+\fr{g(\Gx,v)}{X^+(\Gx,v)}, 
\quad (\Gx,v)\in\CL.
\label{4.38}
\eeq
Next we introduce the singular integral with the kernel $dV$ on the torus $\CR$
\beq
\Psi(\Gz,u)=\fr{\pi-\Ga}{2\pi i}
\int_{l_0'}\left[\fr{\Gz-\Gz_0}{\Gx-\Gz_0}+\fr{u(\Gz)}{u(\Gx)}\fr{\Gx-\Gz_0}{\Gz-\Gz_0}\right]
\fr{d\Gx}{(\Gx-\Gz)X^+(\Gx,u(\Gx))}.
\label{4.39}
\eeq
Note that this function satisfies the symmetry condition $\Psi(\Gz,u)=\ov{\Psi(\Gz_*,u_*)}$, is discontinuous through the contour 
$l_0'$ with the jump $2(\pi-\Ga)/X^+(\Gx,v)$ and continuous through the contours $l_0''$ and $l_1$.

We apply the Sokhotski-Plemelj formulas, the continuity principle and the generalized Liouville
theorem on the torus $\CR$ to deduce the general solution to the Riemann-Hilbert problem
(\ref{4.12})
\beq
\GF(\Gz,u)=X(\Gz,u)[\Psi(\Gz,u)+\GO(\Gz,u)], \quad (\Gz,u)\in\CR.
\label{4.40}
\eeq
Here, $\GO(\Gz,u)$ is a rational function on the torus $\CR$. Its form is determined from the following 
additional conditions the solution has to satisfy.

(i) Since the kernel $dV$  has
simple poles at the points $(\Gz_0,\pm u(\Gz_0))$ on both sheets of the surface $\CR$,
the function $\Psi(\Gz,u)$ has simple poles at these points.
This gives simple poles to the function $\GF(\Gz,u)$ which are unacceptable.
To remove them, we admit that the function $\GO(\Gz,u)$ also has  simple poles at the points 
$(\Gz_0,\pm u(\Gz_0))$ and require that
\beq
 \mathop{\rm res}\limits_{\Gz=\Gz_0}[\Psi(\Gz,u)+\GO(\Gz,u)]=0.
\label{4.41}
\eeq

(ii) Owing to the simple pole of the function $X(\Gz,u)$
at the point $\Bq_1=(\Gz_1,u(\Gz_1))\in\CR$, the function $\Psi(\Gz,u)+\GO(\Gz,u)$ vanishes at this point,
\beq
\Psi(\Gz_1,u(\Gz_1))+\GO(\Gz_1,u(\Gz_1))=0.
\label{4.42}
\eeq

(iii) Because of the simple zero of the function $X(\Gz,u)$
at the point $\Bp_1=(\Gn_1,u(\Gn_1))\in{\Bbb C}_1$,  the function $\GO(\Gz,u)$
has a simple pole at this point.

(iv) Since the functions $X(\Gz,u)$ and $\Psi(\Gz,u)$ are bounded at infinity
and the function $\GF(\Gz,u)$ has to be bounded, the function
$\GO(\Gz,u)$ must be bounded as $\Gz\to\infty$ as well. 

(v) Owing to the symmetry of the functions $\GF(\Gz,u)$, $X(\Gz,u)$, and $\Psi(\Gz,u)$,
the function $\GO(\Gz,u)$ has to be symmetric,  $\GO(\Gz,u)=\ov{\GO(\Gz_*,u_*)}$.

The most general form of such a rational  function is
$$
\GO(\Gz,u)=M_0+(M_1+iM_2)\fr{u(\Gz)+u(\Gn_1)}{\Gz-\Gn_1}-
(M_1-iM_2)\fr{u(\Gz)-u(\bar\Gn_1)}{\Gz-\bar\Gn_1}
$$
\beq
+(M_3+iM_4)\fr{u(\Gz)+u(\Gz_0)}{\Gz-\Gz_0}-
(M_3-iM_4)\fr{u(\Gz)-u(\Gz_0)}{\Gz-\Gz_0},
\label{4.43}
\eeq
where $M_j$ ($j=0,1,\ldots,4$) are arbitrary  real constants.
This function satisfies the conditions (iii) and (v). Because of the growth of the function
$u(\Gz)$ at infinity as $\Gz^{3/2}$, in general, the function $\GO(\Gz,u)$ is unbounded at infinity.
It becomes bounded and the condition (iv) is satisfied if 
\beq
M_2=-M_4.
\label{4.44}
\eeq
To fulfill the property (i), that is to meet the condition (\ref{4.41}),  define the  constants $M_3$ and $M_4$ as 
\beq
M_3=0, \quad M_4=\fr{\pi-\Ga}{4\pi}\int_{l_0'}\fr{d\Gx}{vX^+(\Gx,v)}.
\label{4.45}
\eeq
 Finally, we make sure that the rational function  $\GO(\Gz,u)$ has property (ii).
 Denote
 $$
 R+iJ=\Psi(\Gz_1,u(\Gz_1)), \quad R_0+iJ_0=\fr{2u(\Gz_1)}{\Gz_1-\Gz_0},
 $$
 \beq
 R_1+iJ_1=\fr{u(\Gz_1)+u(\Gn_1)}{\Gz_1-\Gn_1},\quad 
  R_2+iJ_2=\fr{u(\Gz_1)-u(\bar\Gn_1)}{\Gz_1-\bar\Gn_1}.
  \label{4.46}
  \eeq
  where $R$, $J$, $R_l$, and $J_l$ ($l=0,1,2$) are real.
Then equation (\ref{4.42}) implies 
$$
M_1=\fr{(R_0-R_1-R_2)M_2-J}{J_1-J_2},
$$
\beq
M_0=(R_2-R_1)M_1+(J_1+J_2-J_0)M_2-R.
\label{4.47}
\eeq
 Having uniquely obtained all the five constants $M_j$ ($j=0,1,\ldots,4$) in the rational function
 $\GO(\Gz,u)$,  we simplify it as
 \beq
\GO(\Gz,u)=M_0+(M_1+iM_2)\fr{u(\Gz)+u(\Gn_1)}{\Gz-\Gn_1}-
(M_1-iM_2)\fr{u(\Gz)-u(\bar\Gn_1)}{\Gz-\bar\Gn_1}-\fr{2iM_2u(\Gz)}{\Gz-\Gz_0}
\label{4.47'}
\eeq
 and complete the solution procedure for the Riemann-Hilbert problem
 (\ref{4.12}) by analyzing the behavior of the function $\GF(\Gz,u)$ at the points
 $a$ and $\infty$.

Analysis of the integral $\Psi(\Gz,u)$ with the kernel $dV$ 
as $\Gz\to a$ and $\Gz\to \bar a$ yields
 \beq
 \Psi(\Gz,u)\sim \fr{\pi-\Ga}{\Gs\pi i X^+(a,u(a))}\log(\Gz-a), \quad \Gz\to a, 
 \quad  \Psi(\Gz,u)\sim\const, \quad \Gz\to \bar a,
 \label{4.48}
 \eeq
 where $\Gs=1$ when $a\ne m$ and $\Gs=2$ otherwise. Here, we took into account that the starting
 and terminal points of the contour $l_0'$ are $+\infty-i0$ and $a\in l_0^-$, respectively.
 Therefore, in virtue of  (\ref{4.40}), the function $\GF(\Gz,u)$ has the logarithmic singularity,
 as $\Gz\to a$, and its asymptotics coincides with (\ref{4.15})
 predicted by employing formulas (\ref{2.4}) and   (\ref{2.5}) and implementing asymptotic analysis of the conformal map $z=f(\Gz)$ in Section \ref{S3}. As $\Gz\to \bar a$, the function $\GF(\Gz,u)$ is bounded.
 
 Analyze next the behavior of the functions  $\Psi(\Gz,u)$ and  $\GF(\Gz,u)$ at infinity.
 By making the substitutions $\tau=1/\Gx$, $t=1/\Gz$ we deduce as $\Gz\in l_0'=(a,+\infty)^-$,
\beq
\Psi^+\left(
\fr{1}{t},u\left(\fr{1}{t}
\right)\right)
=\fr{\pi-\Ga}{X^+(1/t,u^-(1/t))}
+\fr{\pi-\Ga}{2\pi i}(\CI_1+\CI_2),
\label{4.49}
\eeq
where
$$
\CI_1=\int_{0}^{1/a}
\fr{1-t\Gz_0}{1-\tau\Gz_0}
\fr{d\Gj}{(\Gj-t)X^+(1/\Gj,u^-(1/\Gj))},
$$
\beq
\CI_2=\int_{0}^{1/a}
\sqrt{\fr{t(1-t)(1-mt)}{\Gj(1-\Gj)(1-m\Gj)}}\fr{1-\Gj\Gz_0}{1-t\Gz_0}
\fr{d\Gj}{(\Gj-t)X^+(1/\Gj,u^-(1/\Gj))}.
\label{4.50}
\eeq
As $t\to 0$ and $\Gz\in l_0'$, the integral $\CI_2$ is bounded, whilst the 
integral $\CI_1$ has a logarithmic singularity, and we have
\beq
\Psi^+(\Gz,u)\sim \fr{\pi-\Ga}{2\pi  i X^+(+\infty-i0,u^-)}\log \Gz, \quad \Gz\to +\infty-i0,
\label{4.51}
 \eeq
 where $X^+(+\infty-i0,u^-)=i|(\Gz_1-\Gz_0)/(\Gn_1-\Gz_0)|$. It turns out that the principal term of the asymptotics 
 of  the function $\Psi(\Gz,u)$ as $\Gz\to+\infty+i0$ along the contour $l_0''$ is the same as 
 in (\ref{4.51}).
 In virtue of formula (\ref{4.40}) we discover that $\GF(\Gz,u)\sim (\pi-\Ga)/(2\pi i)\log\Gz$ as 
 $\Gz\to\infty$
 along the contour $l_0$. This formula is identical to the asymptotics (\ref{4.14})
 obtained before by means of the asymptotic analysis of the conformal mapping in Section \ref{S3}.

 \setcounter{equation}{0}
 
\section{Exact representation of the conformal mapping $z=f(\Gz)$}

The conformal mapping from the parametric domain $\CD$ onto the flow domain $\tilde \CD$
is given by
\beq
f(\Gz)=i\int_a^\Gz\fr{N_0+i N_1\Gx}{p^{1/2}(\Gx)} e^{-i\GF^+(\Gx,u(\Gx))}d\Gx,\quad \Gz\in \CD,
\label{5.1}
\eeq
The solution obtained has to satisfy two extra 
equations, the circulation condition (\ref{2.0})
and the complex equation
\beq
\int_{l_1} \fr{df}{d\Gz}d\Gz=0
\label{5.3}
\eeq
which guarantees that the conformal mapping is single-valued.
Owing to formulas (\ref{2.4}) and (\ref{2.5}) we  transform the circulation condition (\ref{2.0}) 
 into the form
 \beq
 \int_0^1\fr{(N_0+N_1\Gx)d\Gx}{\sqrt{|p(\Gx)|}}=-\fr{\GG}{2U}.
 \label{5.4}
 \eeq
 This integral can be computed and expressed though the complete elliptic integrals of the first
 and second kind. We have
 $$
  \int_0^1\fr{d\Gx}{\sqrt{|p(\Gx)|}}=2k\BK(k),
$$
 \beq
 \int_0^1\fr{\Gx d\Gx}{\sqrt{|p(\Gx)|}}=\fr{2}{k}[\BK(k)-\BE(k)], \quad k=\fr{1}{\sqrt{m}}.
 \label{5.5}
 \eeq
 On substituting the integrals (\ref{5.5}) into equation  (\ref{5.4}) we obtain a relation between 
 the two real constants $N_0$ and $N_1$,
 \beq
 N_0=-\fr{1}{k\BK(k)}
 \left[\fr{\GG}{4U}+\fr{\BK(k)-\BE(k)}{k}N_1\right].
 \label{5.6}
 \eeq
 Transform now the second condition, equation (\ref{5.3}). On substituting the representation (\ref{5.1})
 into equation (\ref{5.3}) and denoting
 \beq
 c^\pm(\Gx)=\cos[\GF^+(\Gx, u^\pm(\Gx))], \quad  
 s^\pm(\Gx)=\sin[\GF^+(\Gx, u^\pm(\Gx))],\quad \Gx\in l_1^\pm,
 \label{5.7}
 \eeq
 we write equation (\ref{5.3}) in the form
 $$
 N_0\int_0^1\fr{c^+(\Gx)+c^-(\Gx)-is^+(\Gx)-is^-(\Gx)}{\sqrt{|p(\Gx)|}}d\Gx
 $$
 \beq
 +N_1 \int_0^1\fr{c^+(\Gx)+c^-(\Gx)-is^+(\Gx)-is^-(\Gx)}{\sqrt{|p(\Gx)|}}\Gx d\Gx=0.
 \label{5.8}
 \eeq
Note that since $\I\GF^+(\Gx,u^\pm(\Gx))=0$ on $l_0^\pm$, the functions $c^\pm(\Gx)$
and  $s^\pm(\Gx)$ are real-valued. On taking the real and imaginary parts of equation (\ref{5.8})
we deduce 
\beq
N_0 C_0+N_1C_1=0
\label{5.8'}
\eeq
and 
\beq
 C_0S_1-C_1S_0=0.
\label{5.9}
\eeq 
Here,
\beq
C_j=\int_0^1\fr{c^+(\Gx)+c^-(\Gx)}{\sqrt{|p(\Gx)|}}\Gx^jd\Gx, \quad 
S_j=\int_0^1\fr{s^+(\Gx)+s^-(\Gx)}{\sqrt{|p(\Gx)|}}\Gx^jd\Gx, \quad j=0,1.
\label{5.10}
\eeq
Thus, in addition to equation (\ref{5.6}) we obtained another equation for the constants $N_0$
and $N_1$ and we have
\beq
N_0=-\left[1-\fr{\BK(k)-\BE(k)}{k^2\BK(k)}\fr{C_0}{C_1}\right]^{-1}\fr{\GG}{4Uk\BK(k)},
\quad N_1=-\fr{N_0C_0}{C_1}.
\label{5.11}
\eeq 
Equation (\ref{5.9}) is a real transcendental equation for the two real parameters
$a$ and $m$. This enables us to conclude that the original  problem 
possesses a one-parametric family of solutions provided equation (\ref{5.9})
has a solution and it is unique.

\setcounter{equation}{0}

\section{Reconstruction of the vortex boundary and the wedge sides.  Numerical results}

The first step of the  numerical procedure is to select a value of the parameter $m$
and determine the parameter $a$ by solving the transcendental equation (\ref{5.9}).
At this stage, we need to evaluate the integrals 
(\ref{5.10}) whose integrands depend on the solution $\GF^+(\Gx,u)$ on the contour $l_1$
 and are independent of the constants $N_0$ and $N_1$. We note that we need the solution only on the first sheet.
The function  $\GF^+(\Gx,u)$ is given by (\ref{4.40}) in terms of the factorization function 
$X(\Gz,u)$, the integral $\Psi(\Gz,u)$, and the rational function $\GO(\Gz,u)$.
The factorization function depends on the point $\Bq_1$ and the integer $n_a$,
the solution to the Jacobi inversion problem.
It turns out that for all realistic values of the problem parameters used in the numerical tests and
when $\Bp_1=(\Gn_1,u(\Gn_1))\in {\Bbb C}_1$, the point $\Bq_1=(\Gz_1,u(\Gz_1))$
also lies in the upper sheet ${\Bbb C}_1$, and $n_a=0$. We remark that the final solution, the conformal
map $z=f(\Gz)$, is independent of not only the  point $\Bp_1$ but also the real point $\Gz_0<0$.

Describe first the algorithm for the function $X(\Gz,u)$. It is convenient to write it in the form
\beq
X(\Gz,u)=\Gl\exp\left\{\fr{u(\Gz)}{i(\Gz-\Gz_0)}\left(\fr12-2n_a\right)\CJ_1(\Gz)-\fr{\CJ_2(\Gz,u)+\CJ_3(\Gz,u)}{2}\right\},
\label{6.1}
\eeq
where $\Gl=1$ if $(\Gz,u)\in {\Bbb C}_1\setminus l_0$,  $\Gl=i$ if 
$(\Gz,u)\in l_0\subset{\Bbb C}_1$,
$$
\CJ_1(\Gz)=\int_m^\infty
\fr{(\Gx-\Gz_0)d\Gx}{\sqrt{\Gx(\Gx-1)(\Gx-m)}(\Gx-\Gz)},
$$
\beq
\CJ_2(\Gz,u)=
\int_\Gg\left(\fr{\Gz-\Gz_0}{\Gx-\Gz_0}+\fr{u(\Gz)}{u(\Gx)}\fr{\Gx-\Gz_0}{\Gz-\Gz_0}\right)
\fr{d\Gx}{\Gx-\Gz},
\quad
\CJ_3(\Gz,u)=
\int_\Gg\left(\fr{\Gz-\Gz_0}{\bar\Gx-\Gz_0}-\fr{u(\Gz)}{u(\bar\Gx)}\fr{\bar\Gx-\Gz_0}{\Gz-\Gz_0}\right)\fr{d\bar\Gx}{\bar\Gx-\Gz}.
\label{6.2}
\eeq
Analyze the first integral. If $\Gz\notin(m,+\infty)$,
by making the substitution $\Gx=1/\tau$ we obtain
\beq
\CJ_1(\Gz)=\fr{1}{\sqrt{m}}\int_0^{1/m}\fr{\CF_1(\tau)d\tau}{\sqrt{\tau(1/m-\tau)}},
\quad \CF_1(\tau)=\fr{1-\Gz_0\tau}{\sqrt{1-\tau}(1-\Gj\Gz)}.
\label{6.3}
\eeq
By employing the integration formula with the Gaussian weights $\pi/n$
and the abscissas $x_j=\cos(j-\fr12)\fr{\pi}{n}$ we write approximately
\beq
\CJ_1(\Gz)\approx\fr{\pi}{\sqrt{m}n}\sum_{j=1}^n\CF_1\left(\fr{1+x_j}{2m}\right).
\label{6.4}
\eeq
When $\Gz\in(m,+\infty)$,  the integral $\CJ_1(\Gz)$ is singular, and its Cauchy principal value
is computed by \cite{ant}
\beq
\CJ_1(\Gz)=-\fr{2\pi\sqrt{m}}{\Gz}\sum_{s=1}^\infty d'_s U_{s-1}\left(\fr{2m}{\Gz}-1\right),
\label{6.5}
\eeq
where
\beq
d'_s=\fr{2}{\pi}\int_{-1}^1
\tilde\CF_1\left(\fr{1+\Gj}{2m}\right)
\fr{T_s(\Gj)d\Gj}{\sqrt{1-\Gj^2}}
\approx\fr{2}{n}\sum_{j=1}^n \tilde\CF_1\left(\fr{1+x_j}{2m}\right)\cos\left(j-\fr12\right)\fr{s\pi}{n},
\quad s=1,2,\ldots.
\label{6.6}
\eeq
$\tilde\CF_1(\Gj)=(1-\Gj\Gz_0)(1-\Gj)^{-1/2}$, and
$T_s(\Gj)$ and $U_s(\Gj)$ are the Chebyshev polynomials. Our numerous tests reveal that if
$\I\Gn_1$ is chosen to be positive, then $\I\Gz_1<0$.
It is convenient to split
the contour $\Gg$ lying in the upper sheet into two segments $(\Gn_1,\Gd)$ and $(\Gd,\Gz_1)$, $\Gz_0<\Gd<0$. 
This leads to the following representations of the integrals $\CJ_2$ and $\CJ_3$:
\beq
\CJ_2(\Gz,u)=j^+(\Gz,u ;\Gz_1)-j^+(\Gz,u ;\Gn_1), \quad \CJ_3(\Gz,u)=j^-(\Gz,u ;\bar\Gz_1)-j^-(\Gz,u ;\bar\Gn_1),
\label{6.7}
\eeq
where $j^\pm(\Gz,u;\Gn)$ are the integrals
\beq
j^\pm(\Gz,u;\Gn)=\int_\Gd^\Gn
\left(\fr{\Gz-\Gz_0}{\Gx-\Gz_0}\pm\fr{u(\Gz)}{u(\Gx)}\fr{\Gx-\Gz_0}{\Gz-\Gz_0}\right)
\fr{d\Gx}{\Gx-\Gz}
\label{6.8}
\eeq
evaluated by the Gauss quadrature formula employed for the integral  (\ref{6.3}). 

Consider next the function $\Psi(\Gz,u)$.  If $\Gz\notin l_0'=(a,+\infty)^-$, then the integral is not singular. On making the substitution $\Gx=1/\Gj$, we recast it into
\beq
\Psi(\Gz,u)=-\fr{\pi-\Ga}{2\pi i}\int_0^{1/a}\fr{\CF_*(\Gz,u;\Gj)d\Gj}{\sqrt{\Gj(1/a-\Gj)}},
\label{6.9}
\eeq
where
\beq
\CF_*(\Gz,u;\Gj)=\left(
\fr{(\Gz-\Gz_0)\sqrt{\Gj}}{1-\Gj\Gz_0}-\fr{u(\Gz)(1-\Gz_0\Gj)}{i\sqrt{(1-\Gj)(1-m\Gj)}(\Gz-\Gz_0)}
\right)
\fr{\sqrt{1/a-\Gj}}{(1-\Gj\Gz)X^+(1/\Gj,v^-)},
\label{6.10}
\eeq
where $v^-=u(\Gx-i0)$. This integral is efficiently evaluated by the Gauss quadrature formulas used
in (\ref{6.3}), where $m$ and $\CF_1$ need to be replaced by $a$ and $\CF_*$, respectively.

In the case $\Gz\in l_0'=(a,+\infty)^-$ on the upper sheet, the integral $\Psi(\Gz,u)$ is singular.
We denote $u^-=u(\Gz)$, $\Gz\in l_0'$, and 
by the Sokhotski-Plemelj formulas,
\beq
\Psi^+(\Gz,u^-)=\fr{\pi-\Ga}{X^+(\Gz,u^-)}+\Psi(\Gz,u^-), \quad \Gz\in l_0'.
\label{6.11}
\eeq
We again use the formula \cite{ant} for the Cauchy principal value  $\Psi(\Gz,u^-)$
and similar to the integral $\CJ_1$ obtain
\beq
\Psi(\Gz,u^-)=\fr{\pi-\Ga}{2\pi i}\int_0^{1/a}\fr{\CF(\Gj)d\Gj}{\sqrt{\Gj(1/a-\Gj)}(\Gj-t)}=
-i(\pi-\Ga)a\sum_{s=1}^\infty d_s
U_{s-1}\left(\fr{2a}{\Gz}-1\right), 
\label{6.12}
\eeq
where $t=1/\Gz$, $\Gz=l_0'$,
$$
\CF(\Gj)=
\left[
\fr{(1-\Gz_0t)\sqrt{\Gj}}{1-\Gz_0\Gj}+\sqrt{\fr
{t(1-t)(1-mt)}{(1-\Gj)(1-m t)}}
\fr{1-\Gz_0\Gj}{1-\Gz_0t}
\right]
\fr{\sqrt{1/a-\Gj}}{X^+(1/\Gj,v^-)}.
$$
\beq
d_s
\approx\fr{2}{n}\sum_{j=1}^n \CF\left(\fr{1+x_j}{2a}\right)\cos\left(j-\fr12\right)\fr{s\pi}{n},
\quad s=1,2,\ldots.
\label{6.13}
\eeq
Note that according to (\ref{4.51}) the series in (\ref{6.12}) has a logarithmic singularity as $\Gz\to\infty$.

The third function in the representation (\ref{4.40}) of the solution of the Riemann-Hilbert problem,
the function $\GO(\Gz,u)$,  has three real constants $M_0$, $M_1$, and $M_2$.
The constant $M_2$
is the integral
\beq
M_2=- \fr{\pi-\Ga}{4\pi}\int_{l_0'}\fr{d\Gx}{v^-X^+(\Gx,v^-)}=-\fr{\pi-\Ga}{4\pi i}\int_0^{1/a}\fr{\CF_0(\Gj)d\Gj}
{\sqrt{\Gj(1/a-\Gj)}},
\label{6.14}
\eeq
where
\beq
\CF_0(\Gj)=\fr{\sqrt{1/a-\Gj}}{\sqrt{(1-\Gj)(1-m\Gj)}X^+(1/\Gj,v^-)}.
\label{6.15}
\eeq
This integral is evaluated numerically by the Gauss quadrature rule employed before. The other two constants are expressed through $M_2$ by (\ref{4.47}).

Having equipped with formulas for $\GF^+(\Gx,u)$ efficient for numerical purposes,
we can evaluate the integrals (\ref{5.10}) by applying  the Gauss quadratures
and solve the transcendental equation  (\ref{5.9}). It turns out that, for all numerical tests implemented,
$a=m$. In the following formulas we do not distinguish the points $a$ and $m$ anymore.

\begin{figure}[t]
\centerline{
\scalebox{0.6}{\includegraphics{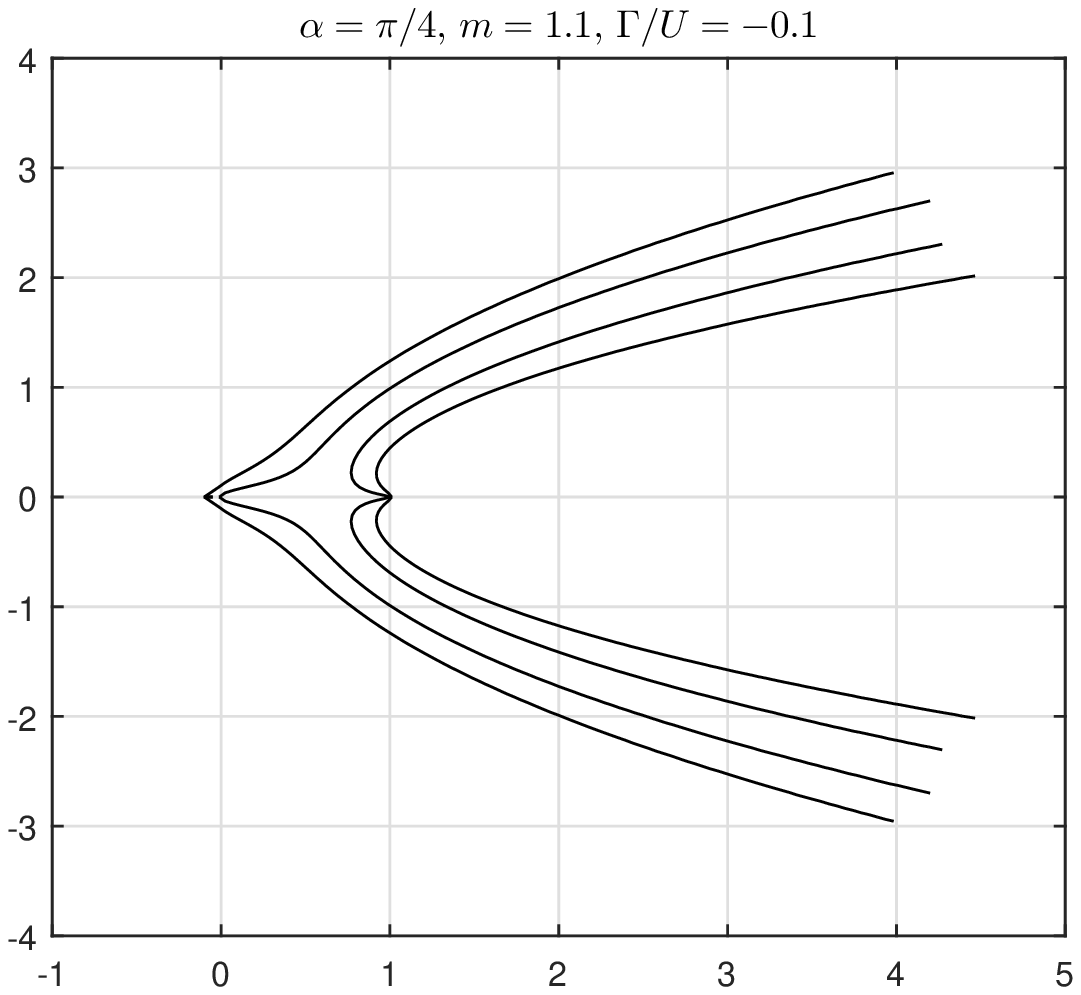}}}
\caption{The preimages of some streamlines when  $\Ga=\pi/4$, $\GG/U=-0.1$, $m=1.1.$}
\label{fig3}
\end{figure}

\begin{figure}[t]
\centerline{
\scalebox{0.5}{\includegraphics{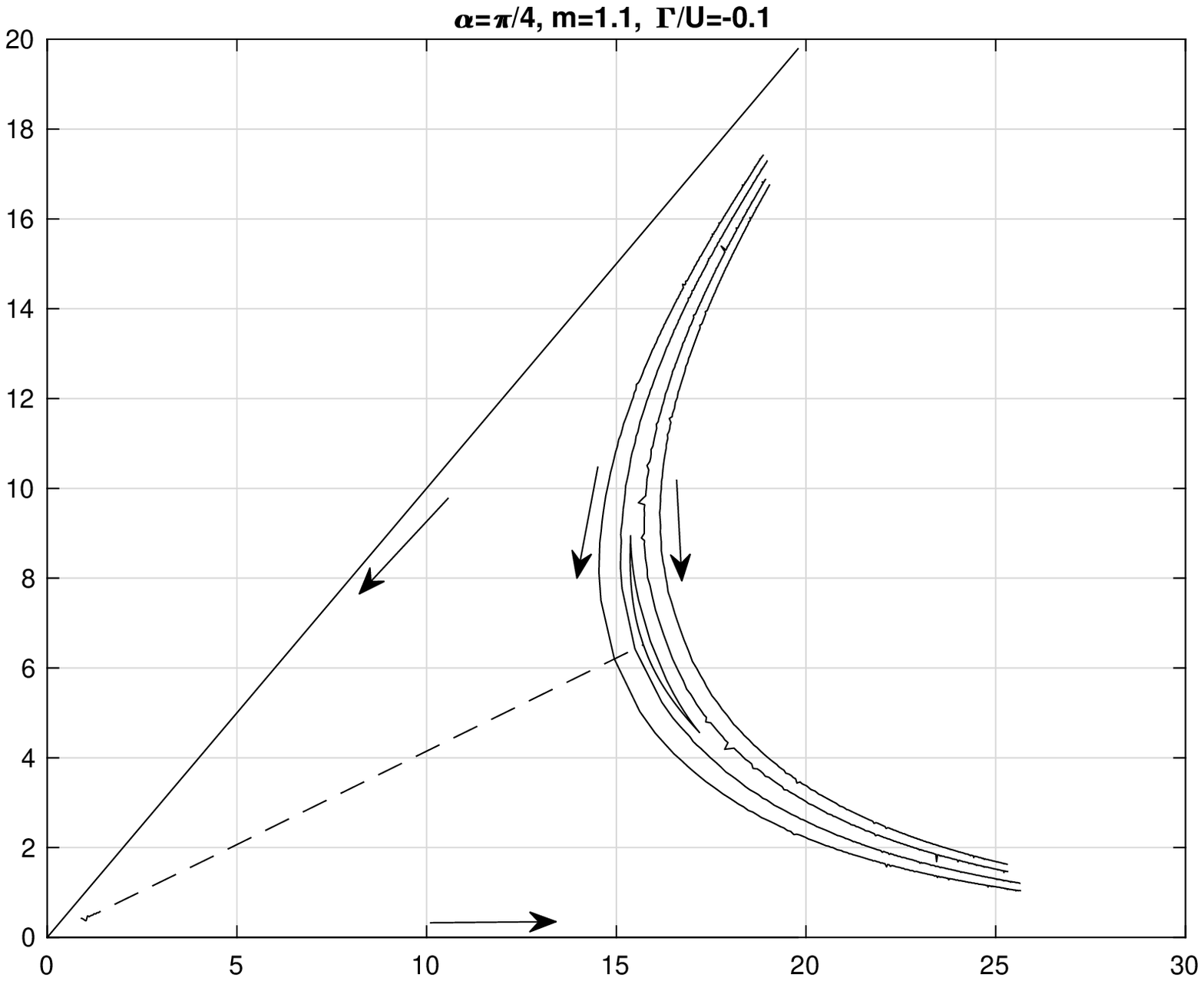}}}
\caption{The vortex domain and some streamlines when  $\Ga=\pi/4$, $\GG/U=-0.1$, $m=1.1.$}
\label{fig4}
\end{figure}

\begin{figure}[t]
\centerline{
\scalebox{0.6}{\includegraphics{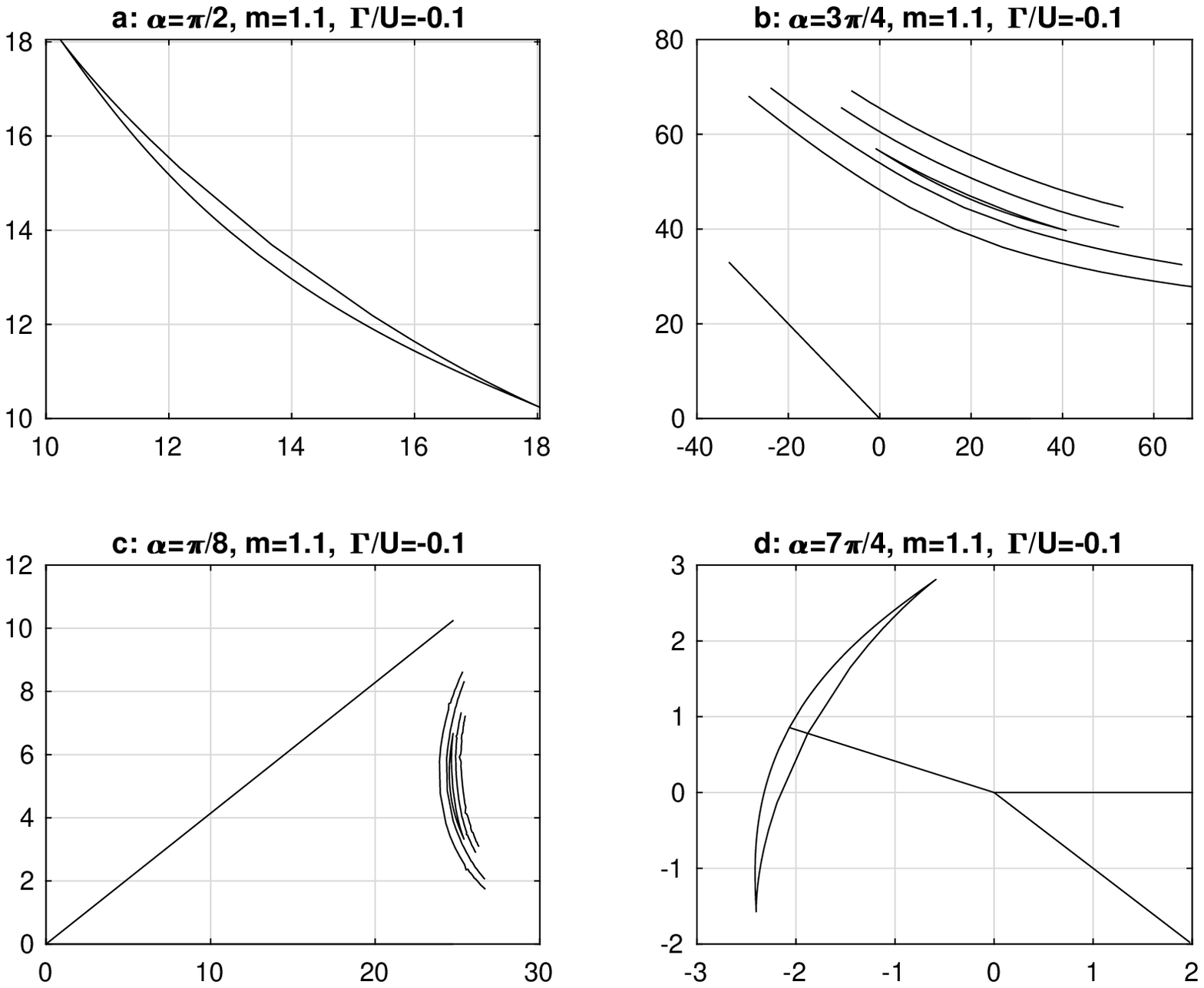}}}
\caption{The vortex domain, the wedge, and some streamlines for $m=1.1$  when
(a) $\Ga=\fr12\pi$, (b)  $\Ga=\fr34\pi$, (c) $\Ga=\fr18\pi$, (d) $\Ga=\fr74\pi$: owing
to penetration of streamlines into the vortex domain the solution does not exist for $\Ga>\pi$.} 
\label{fig5}
\end{figure}

\begin{figure}[t]
\centerline{
\scalebox{0.6}{\includegraphics{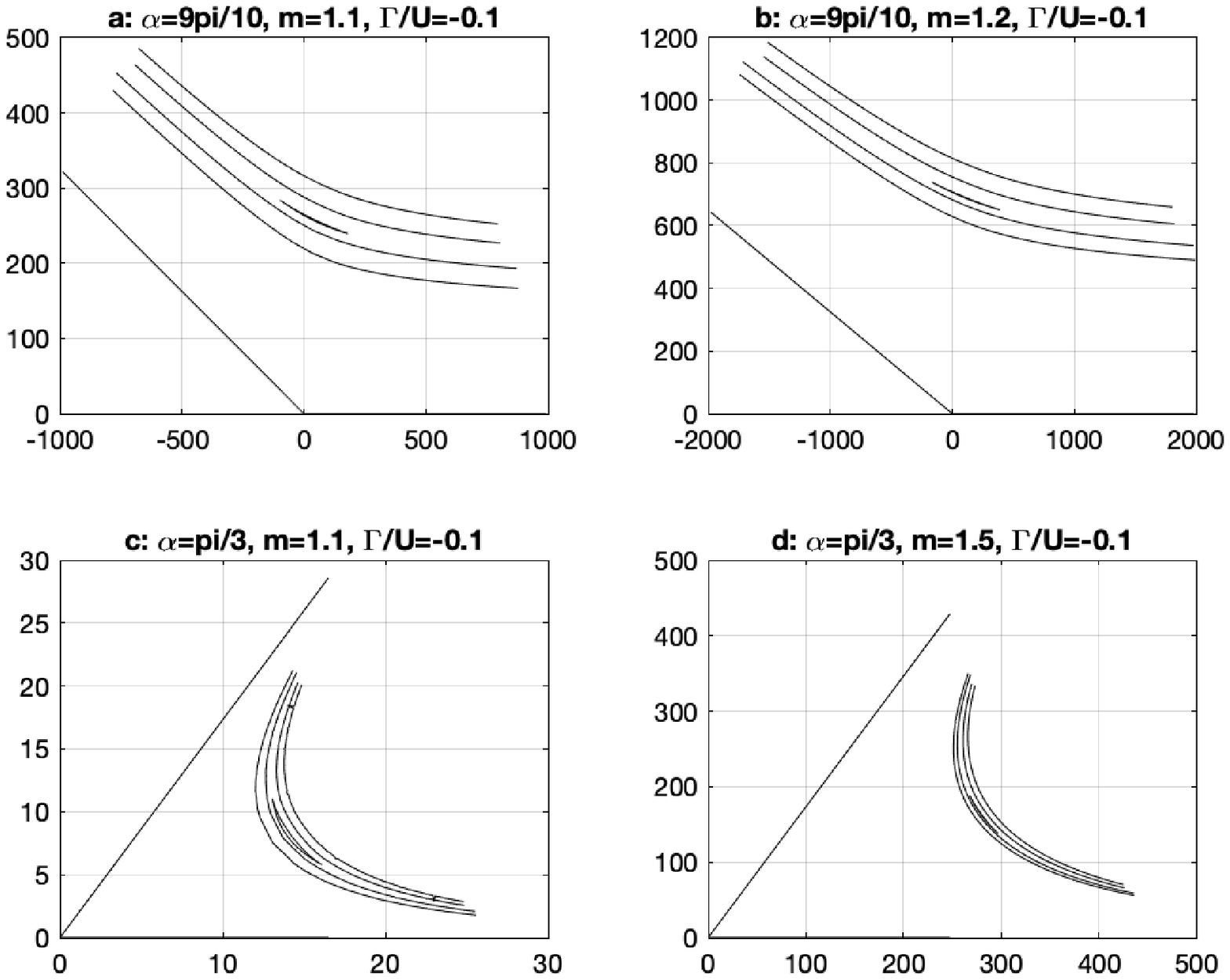}}}
\caption{The vortex domain, the wedge, and some streamlines when
(a) $\Ga=\fr{9}{10}\pi$, $m=1.1$ (b) $\Ga=\fr{9}{10}\pi$, $m=1.2$  (c) $\Ga=\fr13\pi$,  $m=1.1$,
d)  $\Ga=\fr{1}{3}\pi$, $m=1.5$.} 
\label{fig6}
\end{figure} 

 \begin{figure}[t]
\centerline{
\scalebox{0.6}{\includegraphics{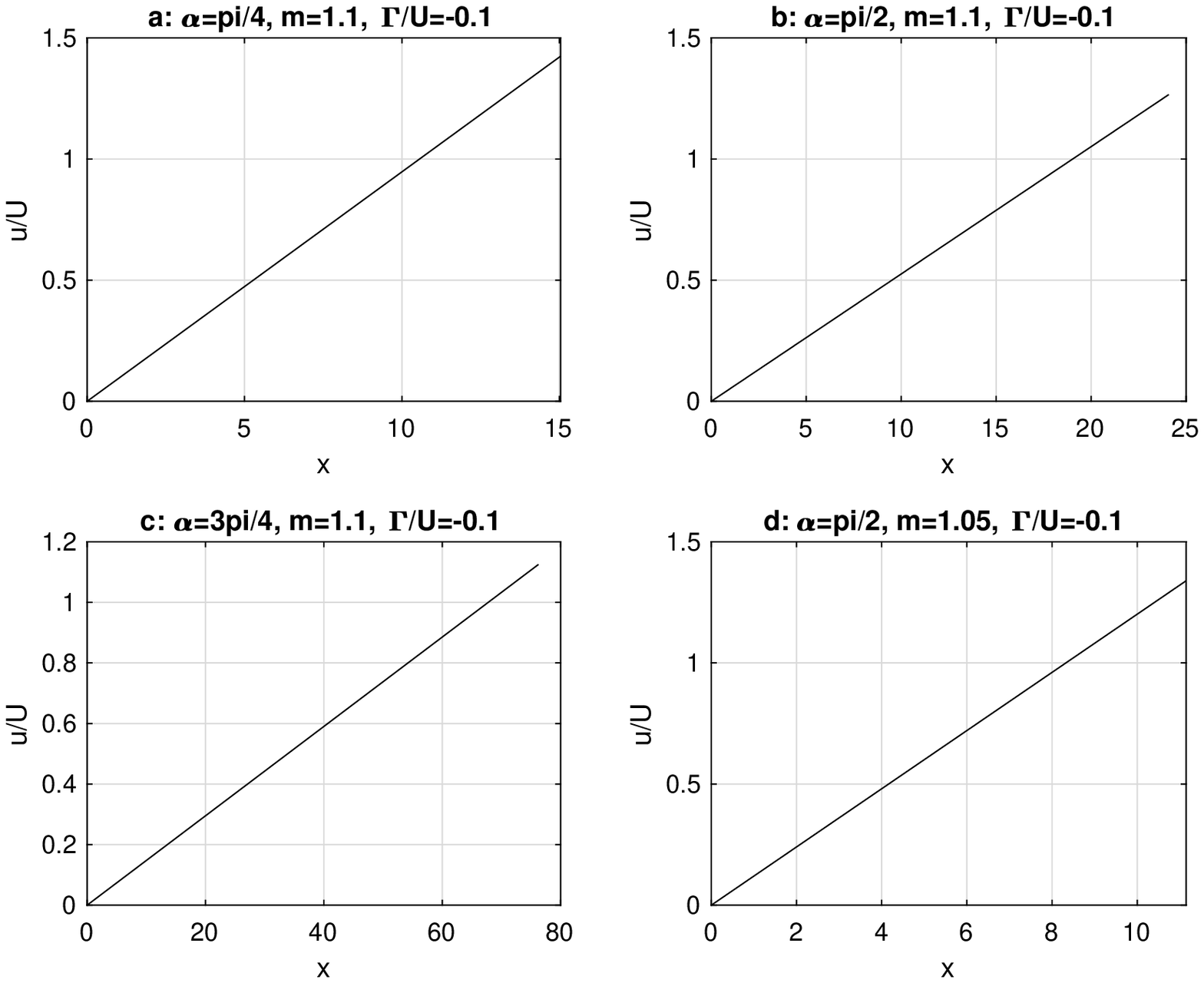}}}
\caption{The normalized tangential velocity $u_x/U$ on the horizontal side of the wedge when $\GG/U=-0.1$
(a) $\Ga=\fr{1}{4}\pi$, $m=1.1$ (b) $\Ga=\fr12\pi$, $m=1.1$  (c) $\Ga=\fr34\pi$,  $m=1.1$,
d)  $\Ga=\fr12\pi$, $m=1.05$.} 
\label{fig7}
\end{figure}

We aim now to show that the conformal map we found indeed
maps the contours $l_0^-$ and $l_0^+$ into the sides
$\arg z=\Ga$ and $\arg z=0$ of the wedge, respectively.
Let $\Gz\in l_0^\pm$. Then
\beq
f(\Gz)=\pm\int_m^\Gz\fr{N_0+N_1\Gx}{\sqrt{|p(\Gx)|}} e^{-i\GF^+(\Gx,v^\pm)}d\Gx\quad  v^\pm=u(\Gx^\pm).
\label{6.16}
\eeq
Notice that on  both sides $l_0^\pm$ of the contour $l_0$, $\R X^+(\Gx,v^\pm)=0$,
$\I\GO(\Gx,v^\pm)=0$. Owing to formula  (\ref{6.11}),
we discover
\beq
 \GF^+(\Gx,v^-)=\pi-\Ga+iP(\Gx,v^-), \quad \Gx\in l_0^-,\quad  
 \GF^+(\Gx,v^+)=iP(\Gx,v^+), \quad \Gx\in l_0^+,
 \label{6.17}
 \eeq
 Here, $P(\Gx,v^\pm)$ are  real-valued  functions,
 \beq
P(\Gx,v^\pm)=\I X^+(\Gx,v^\pm)[\Psi(\Gx,v^\pm)+\GO(\Gx,v^\pm)], \quad \Gx\in l_0^\pm.
\label{6.19}
\eeq
Returning to equation (\ref{6.16}), we deduce
\beq
f(\Gz)=e^{i\Ga}\int_m^\Gz 
\fr{N_0+N_1\Gx}{\sqrt{|p(\Gx)|}}
e^{P(\Gx, v^-)}d\Gx, \quad \Gz\in l_0^-,
\quad
f(\Gz)=\int_m^\Gz \fr{N_0+N_1\Gx}{\sqrt{|p(\Gx)|}} e^{P(\Gx, v^+)}d\Gx, \quad \Gz\in l_0^+.
\label{6.21}
\eeq
 Our numerical tests show that  the function 
 $N_0+N_1\Gx$ is positive for all $\Gx\in(m,+\infty)$. This implies
 that $z=f(\Gz)$ maps the contours $l_0^-$ and $l_0^+$ into the  sides $AB$  and $AC$ of the wedge,
 respectively.

To reconstruct  the vortex boundary, we shall determine the image of the straight path
joining 
the starting point $\Gz=m$ and the terminal point $\Gz=1$.
Analysis of the factorization function $X(\Gz,u)$ and $\Psi(\Gz,u)$ given by (\ref{6.1}) and  (\ref{4.39})
shows that both functions have singularities as $\Gz\to m^-$. To determine the position of the image
of the point $\Gz=1$ with a good accuracy is crucial for the algorithm capability to reconstruct 
the vortex profile and streamlines close to the vortex. Consider the singular integral $\CJ_1(\Gz)$ 
that brings the singularity to the function $X(\Gz,u)$.  Denote $t=1/\Gz$ and for  $1/m<t<1$ we have 
\beq
\CJ_1(\Gz)=\CJ_1^*(t)+\CJ_1^{**}(t), 
\label{6.21.1}
\eeq
where
$$
\CJ_1^*(t)=-t\int_0^{1/m}\left(\fr{1-\Gz_0\Gj}{\sqrt{1-\Gj}}-\fr{1-\Gz_0/m}{\sqrt{1-1/m}}\right)\fr{d\Gj}{\sqrt{\Gj(1-m\Gj)}(\Gj-t)},  \quad \fr{1}{m}<t<1.
$$
\beq
\CJ_1^{**}(t)=-t\fr{1-\Gz_0/m}{\sqrt{1-1/m}}\int_0^{1/m}\fr{d\Gj}{\sqrt{\Gj(1-m\Gj)}(\Gj-t)}, \quad 
\fr{1}{m}<t<1.
\label{6.21.2}
\eeq
The first integral is not singular, and we can manage with it using the Gauss quadrature rule, whilst the second integral that carries the singularity is evaluated exactly,
\beq
\CJ_1^{**}(t)=\fr{1-\Gz_0/m}{\sqrt{1-1/m}}\fr{\pi t}{\sqrt{t(mt-1)}}, \quad \fr{1}{m}<t<1.
\label{6.21.3}
\eeq
By passing to the limit $\Gz\to m^-$ and $\Gz\to m^+\pm i0$ in (\ref{6.1}) and using formula (\ref{6.21.3})
we obtain
\beq
\lim_{\Gx\to m^+\pm i0}X^+(\Gx,v)=\lim_{\Gz\to m^-}X(\Gz,u)=X(m,u)=i\left|
\fr{(\Gz_1-\Gz_0)(\Gn_1-m)}{(\Gz_1-m)(\Gn_1-\Gz_0)}
\right|.
\label{6.21.4}
\eeq

Consider next the function $\Psi(\Gz,u)$ with $(\Gz,u)\in{\Bbb C}_1$.  For $\Gz=1/t\in(1,m)$ we represent it in the form
\beq
\Psi(\Gz,u)=\Psi_1(t)+\Psi_2(t),
\label{6.21.5}
\eeq
where
$$
\Psi_1(t)=\Psi_1^*(t)+\fr{\pi-\Ga}{2\pi i X(m,u)}\log\fr{t-1/m}{t},
$$
\beq
\Psi_2(t)=-\fr{(\pi-\Ga)\sqrt{1-t}}{2(1-\Gz_0 t)} \left[
\fr{\sqrt{t(1-mt)}}{\pi}
\Psi_2^*(t)-r_0 \right], \quad r_0=\fr{1-\Gz_0/m}{\sqrt{1-1/m}X(m,u)}.
\label{6.21.6}
\eeq
The integrals $\Psi_1^*(t)$ and $\Psi^*_2(t)$ are given by
$$
\Psi_1^*(t)=\fr{\pi-\Ga}{2\pi i}\int_0^{1/m}\left[\fr{1-\Gz_0 t}{(1-\Gz_0\Gj)X^+(1/\Gj,v)}-\fr{1}{X(m,u)}\right]\fr{d\Gj}{\Gj-t},
$$
\beq
\Psi_2^*(t)=\int_0^{1/m}\left[\fr{1-\Gz_0 \Gj}{\sqrt{1-\Gj}X^+(1/\Gj,v)}-r_0\right].
\fr{d\Gj}{\sqrt{\Gj(1-m\Gj)}(\Gj-t)}.
\label{6.21.7}
\eeq
They are evaluated numerically in the same manner as before since they are not singular.
On combining these formulas we obtain for $\Gz<m$ close to $m$
\beq
\Psi(\Gz,u)\sim 
\fr{\pi-\Ga}{2\pi i X(m,u)}
\log(m-\Gz)+\fr{\GL_0}{iX(m,u)}, \quad \Gz\to m^-,
\label{6.21.8}
\eeq
where
\beq
\GL_0=\fr{\pi-\Ga}{2}
\left(i-\fr{\log m}{\pi}\right)
+i X(m,u)\Psi_1^*(1/m).
\label{6.21.9}
\eeq
From here we discover the asymptotics of the solution  $\GF(\Gz,u)$ to the Riemann-Hilbert problem (\ref{4.12})
 on the first sheet ${\Bbb C}_1$
and of the corresponding function $\exp\{-i\GF(\Gz,u)\}$ 
$$
\GF(\Gz,u)\sim \fr{\pi-\Ga}{2\pi i}\log(m-\Gz)-i\GL_0+X(m,u)\GO(m,u),
$$
\beq
e^{-i\GF(\Gz,u)}\sim(m-\Gz)^{-1/2+\Ga/(2\pi)}\GL_1, \quad \Gz\to m^-.
\label{6.21.10}
\eeq
Here,
\beq
\GL_1=e^{-\GL_0-iX(m,u)\GO(m,u)}.
\label{6.21.11}
\eeq
The images $f(\Gz)$ of the points $\Gz\in [1,m]$ can now be evaluated with a good accuracy by applying
the formula
\beq
f(\Gz)=i\int_\Gz^m \fr{g(\Gx)d\Gx}{\sqrt{(m-\Gx)(\Gx-\Gz)}}+\fr{2\pi i N}{\Ga}(m-\Gz)^{\Ga/(2\pi)},\quad 1\le\Gz\le m,
\label{6.21.12}
\eeq
where
\beq
g(\Gx)=\fr{(N_0+N_1\Gx)e^{-i\GF(\Gx,u)}}
{\sqrt{\Gx(\Gx-1)}}-N(m-\Gx)^{-1/2+\Ga/(2\pi)},
\quad
N=\fr{(N_0+N_1 m)\GL_1}{\sqrt{m(m-1)}}.
\label{6.21.13}
\eeq
Since the function $g(\Gx)$ is not singular at the point $m$, $g(\Gx)=O((m-\Gx)^{1/2+\Ga/(2\pi)})$, 
$\Gx\to m$, the integral in (\ref{6.21.12}) is evaluated by the Gauss quadrature formula.

At the final step of the procedure we determine  the profile $z=f(\Gz)$ of the vortex domain by employing the 
formula
\beq
f(\Gz^\pm)=f(1)\mp\int_1^{\Gz^\pm}\fr{(N_0+N_1\Gx)e^{-i\GF^+(\Gx^\pm,v^\pm)}d\Gx}{\sqrt{\Gx(1-\Gx)(m-\Gx)}},
 \quad \Gx^\pm, \Gz^\pm\in l_1^\pm,
 \label{6.22}
 \eeq
and the velocity on the wedge walls
\beq
\fr{dw}{Udz}=\fr{u_x-iu_y}{U}=\Gl^\pm\exp\left\{X^+(\Gx^\pm,v^\pm)[\Psi(\Gx^\pm, v^\pm)+\GO(\Gx^\pm,v^\pm)]\right\},
\quad  \Gx^\pm\in l_0^\pm,
\label{6.23}
\eeq
where $v^\pm=u(\Gx^\pm)$, $\Gl^+=1$ as $\Gx\in l_0^+$ ($z\in AC$) and $\Gl^-=e^{i(\pi-\Ga)}$ as $\Gx\in l_0^-$ ($z\in AB$).
 The complex constant $f(1)$ is computed from (\ref{6.21.12}), and the nonsingular integral (\ref{6.22}) is evaluated 
 by the Gauss formula in a simple manner. It is directly verified that the argument  of the exponential function
 in (\ref{6.23}) is real, and the boundary conditions (\ref{2.1}) are satisfied. Analysis of the asymptotics of the functions $\Psi(\Gx^\pm,v^\pm)$ as $\Gx^\pm\to m$ yields
 \beq
 \Psi(\Gx^\pm,v^\pm)\sim \fr{\pi-\Ga}{2\pi i X(m,u)}\log(\Gz-m), \quad \Gx^\pm\to m,
 \label{6.24}
 \eeq
 and therefore 
\beq
\fr{dw}{Udz}\sim 
K(\Gx-m)^{1/2-\Ga/(2\pi)}, \quad \Gx\to m, \quad \Gx\in l_0^\pm,  \quad K=\const.
\label{6.25}
\eeq
This implies $dw/dz\sim cz^{\pi/\Ga-1}$, $z\to 0$, that is consistent with the analysis of Section 4.1  in the 
case $a=m$.

For numerical tests, we select
 $ \Gz_0=-5$,  $\Gd=-2$, and $\Gn_1=1+i/2$. We emphasize that all the curves presented in Figures 3  to  7
 are independent of the choice of these parameters provided $-\infty<\Gz_0<\Gd<0$ and $\I\Gn_1\ne 0$.
 In the case $m=1.1$,
the output of the Jacobi inversion problem is 
$n_a=n_b=0$ and $\Gz_1=0.668204-0.835295 i\in {\Bbb C_1}$. 
For all other values of the 
parameters we tested, the situation is similar: $n_a=n_b=0$ and  $\Gz_1 \in {\Bbb C_1}$ with $\I\Gz_1<0$.

The only one restriction for the parameter $\GG/U$ is it has to be negative. If it is positive,
then the model problem  (\ref{2.0}) to (\ref{2.3}) is solvable if the direction of the velocity
on the wedge walls is reversed and  the boundary condition (\ref{2.1})
is changed accordingly. Otherwise the solution does not
exist. The parameter $\GG/U$ is selected to be $-0.1$ for all our tests.
The preimages 
\beq
\I\int_{\Gx_0}^\Gz \Go_0(\Gn)d\Gn=0
\label{6.26}
\eeq
in the $\Gz$-parametric plane   of some streamlines around the vortex are presented
in Fig.3. Here, $\Go_0(\Gn)$ is given by (\ref{3.2}), $\Gx_0<0$ for the preimages of streamlines facing the 
infinite point of the wedge
and $1<\Gx_0<m$ of those facing the corner.
The curves are symmetric with respect to the real axis. In Fig.3, the point $\Gx_0$ is chosen to be
$-0.01$, $-0.05$ and $1.001$, $1.01$, respectively, while $\alpha=\fr14\pi$ and $m=1.1$.
The vortex and the streamlines corresponding to their preimages in Fig. 3 for the same $\Ga$ and $m$ are shown in Fig.4.
Two cusps are seen at the points where the streamline branches and where the branches merge.
 The dashed line demonstrates the preimage of the straight path from $m$ to $1$.
Figures 5 a-b show the vortex domains and in cases b and c some streamlines around the vortex when $m=1.1$ and the angle
$\Ga=\fr12\pi$ (a), $\Ga=\fr34\pi$ (b), and $\Ga=\fr18\pi$ (c). It turns out that for angles $\Ga>\pi$
the preimage of the straight path
from $\Gz=m$ to $\Gz=1$ intersects the preimage of the contour $l_1$ (Fig. 5 d). It is also found
that the streamlines close to the vortex boundary penetrate it. These results indicate that
a vortex domain with the prescribed properties does not exist in wedges when $\Ga>\pi.$
Some other samples of vortices and streamlines nearby are given in Fig. 6 a-d. 
In Fig. 6 a and b, $\Ga=\fr{9}{10}\pi$ with $m=1.1$ and $m=1.2$.  In Fig. 6 c and d,  $\Ga=\fr13\pi$
with $m=1.1$ and $m=1.5$. It is seen that when the conformal mapping free parameter $m$ grows
the vorticity  drifts out of the wedge corner. Fig. 6 a and b also show that when $\Ga$ grows and approaches $\pi$, 
the vortex sheet length grows as well,  and the vortex domain moves away from the corner. 
The normalized tangential velocity $u_x/U$ on the horizontal side of the wedge is plotted for some $\Ga$ and $m$ in Fig.7.
It is seen that when the angle $\Ga$ decreases, the speed $u_x$ is increasing. It is also observed that when $m$
decreases and the vortex becomes closer to the corner, the speed grows.

 \vspace{.1in}

\setcounter{equation}{0}

\section{Conclusions} 

A two-dimensional flow of an inviscid incompressible fluid in a wedge with impenetrable walls  about a Sadovskii vortex  has been investigated by the method of conformal mappings. The shape of the vortex domain is unknown {\it a priori},
and no stagnation points on the wedge sides are admitted by the model. The map has been constructed by solving a symmetric  Riemann-Hilbert problem on a finite and a semi-infinite segments lying on an elliptic Riemann surface.
Owing to the asymptotics at infinity  of the classical Weierstrass analogue of the Cauchy kernel for an elliptic surface  and the right hand-side of the Riemann-Hilbert problem, the use of the Weierstrass kernel leads to divergent integrals.
A new analogue of the Cauchy kernel on a hyperelliptic surface with properties required has been proposed, and a closed-form solution  to the 
Riemann-Hilbert problem has been constructed.  It is noted that the Weierstrass kernel would have worked
 if another map from the exterior of two finite segments on the flow domain had been applied. 
 The initial authors' efforts  to numerically implement the procedure and recover such a map failed.
On the contrary, the mapping from a finite  and a semi-infinite two-sided segments onto the vortex
boundary and the wedge walls, respectively, has been numerically successful.  
 The conformal mapping recovered by solving the Riemann-Hilbert problem  has a free geometric
 parameter and
 two problem parameters, the wedge angle $\Ga$ and the ratio $\GG/U$ of the circulation and
 the speed on the vortex boundary. The  vortex  and the wedge boundaries have been determined by quadratures and 
numerically reconstructed. The vortex boundary formed by two branches of the same streamline 
in all numerical tests share the same feature of having
two cusps. 
It has also been found that although
for the wedge angle $\Ga$ greater than $\pi$, it is possible to reconstruct a vortex domain by solving  the free boundary
problem, the solution cannot be accepted as a physical one since streamlines close to the vortex domain 
penetrate he vortex interior. Thus the solution exists and has one free real parameter for $\Ga<\pi$ and  does not exist
when $\Ga>\pi$.

\vspace{.1in}

{\bf Acknowledgments.} 

The authors thank
the Isaac Newton Institute for Mathematical Sciences, Cambridge, for support and hospitality during the program Complex analysis: techniques, applications and computations,  where a part of work on this paper was undertaken.

\end{document}